\documentclass[preprint,11pt,authoryear]{elsarticle}

\usepackage{setspace}
\usepackage[margin=3cm]{geometry}

\usepackage{url}

\usepackage{tabularx}
\usepackage{xltabular}

\usepackage{algorithm}
\usepackage{algpseudocode}

\usepackage{tikz}

\usepackage{amssymb}
\usepackage{amsmath}

\journal{European Journal of Operational Research}

\begin{document}

\begin{frontmatter}



\title{On an optimization model for firefighting helicopter planning} 

	\author[label1,label3]{Marta Rodríguez Barreiro}
\author[label2,label3]{María José Ginzo Villamayor} 
\author[label4]{Fernando Pérez Porras}
\author[label1]{María Luisa Carpente Rodríguez} 
\author[label1,label5]{Silvia María Lorenzo Freire}   

\affiliation[label1]{organization={MODES Group, Department of Mathematics, University of A Coruña},
	country={Spain}}
\affiliation[label2]{organization={Modestya Group, Department of Statistics, Mathematical Analysis and Optimization, University of Santiago de Compostela},
	country={Spain}}
\affiliation[label3]{organization={Galician Centre of Mathematical Research and Technology (CITMAga)},
	country={Spain}}
\affiliation[label4]{organization={Department of Graphic and Geomatic Engineering, University of Córdoba},
	country={Spain}}
\affiliation[label5]{organization={Centre for Information and Communications Technology Research (CITIC)},
	country={Spain}}


\begin{abstract}
			During a wildfire, the work of the aerial coordinator is crucial for the control of the wildfire and the minimization of the burned area and the damage caused. Since it could be very useful for the coordinator to have decision-making tools at his/her disposal, this framework deals with an optimization model to obtain the optimal planning of firefighting helicopters, deciding the points where the aircraft should load water, the areas of the wildfire where they should work, and the rest bases to which each helicopter should be assigned. It was developed a Mixed Integer Linear Programming model which takes into account the configuration of helicopters, in closed circuits, as well as the flight aerial regulations in Spain. Due to the complexity of the model, two algorithms are developed, based on the Simulated Annealing and Iterated Local Search metaheuristic techniques. Both algorithms are tested with real data instances, obtaining very promising results for future application in the planning of aircraft throughout a wildfire evolution.
\end{abstract}


\end{frontmatter}



\section{Introduction}
\label{introduction}
Wildfires are one of the greatest problems facing humanity nowadays. The European Forest Fire Information System (EFFIS) reported than more than 467000 hectares were burned in Europe in 2023. Climate change, with the resulting increase in temperatures and extreme weather conditions, has led to wildfires consuming more and more hectares, resulting in what are known as large wildfires (those measuring more than 500 hectares). A report by the World Wildlife Fund (WWF; \cite{WWF}) warns that the new reality of 21st century wildfires is that they are fewer in number, but larger. Furthermore, due to drought conditions caused by climate change, the wildfire season is starting earlier. 

In the study developed by \cite{INFOCA} about the management of aerial firefighting resources in Andalucía (Spain), the authors state that many aerial resources are working at the same time to extinguish a wildfire, with the simultaneous action of more than 10 aerial resources being frequent. The aerial coordinator is in charge of managing the operations of all these aircraft and, according to these authors, it is more effective when a global management is made for the whole wildfire, thus increasing the safety and efficiency of all operations. The aerial coordinator must make decisions about the extinguishing aircraft, their working periods, working areas on the wildfire, water loading points, etc. All of these decisions must be made quickly, once the wildfire is detected, and in the context of an emergency situation with all the associated stress. 

Much effort is being put into developing tools to assist in the planning of aircraft firefighting operations due to the complexity and severity of the problem. In \cite{Granda} a review of recent literature related to wildfire suppression is presented. Following a classification previously established in \cite{Martell}, they divide the wildfire suppression process into four phases: resource acquisition and strategic deployment, resource mobilization, initial fire attack and extended attack management. 

Table \ref{tab:review_characteristics} summarizes the characteristics of some works found in the literature. ``Phase'' refers to the previously mentioned classification, where 1 corresponds to resource acquisition and strategic deployment, 2 to resource mobilization, 3 refers to initial fire attack, and 4 to extended attack management. ``Stochastic'' refers to whether it is a stochastic model, ``Graph'' indicates if the model is represented with a graph, and ``Grid'' refers to the terrain representation. 
This terrain representation consider cells in a grid to represent the different portions of terrain. ``Fireline'' refers to the concept of fireline production. In these models, extinguishing resources are associated with a fireline production, and the wildfire is contained when it is completely enclosed by this fireline. This type of modeling may not entirely conform to reality, particularly when working with a single type of resource. In the case of this work, where terrestrial resources are not considered, this concept cannot be used since the perimeter of the wildfire may not be enclosed by the perimeter formed by the aerial extinguishing resources, because they cannot work in the areas where the terrestrial means are located. 

``Dynamic'' refers to whether the model takes into account the dynamic nature of the wildfire, as it varies in both location and intensity over time. Lastly, ``Legislation'' refers to any aspect related to the work of the resources, such as maximum working times.

\begin{table}[h]
	\centering
	\footnotesize
	\begin{tabularx}{\textwidth}{cccccccc}
		Reference & Phase & Stochastic  & Graph & Grid	& Fireline & Dynamic & Legislation \\
		\hline
		\rule{0pt}{3ex}    \cite{Mendes_Alvelos} & 3 & & x & & & x & \\
		\cite{Mendes_Alvelos_2} & 3 & & x & & & x & \\ 
		\cite{Belval} & 3 & x &  & x & & x & \\
		\cite{Wiitala} & 3 & & & & x & x & \\
		\cite{Ntaimo} & 1, 3 & x & & & x & x & \\
		\cite{Donovan} & 4 & & & & x & x & \\
		\cite{Zhou} & 4 & x & & x &  & x & \\
		\cite{Jorge_1} & 4 & & & & x & x &  x\\
		\cite{Jorge_descomposicion} & 4 & & & & x & x & x\\
		\cite{shahidi} & 3 & & x & & &  & \\
		\cite{Skorin} & 4 & &  & & & x & x\\
		Present work & 4 &  & x & &  & x & x \\
		\hline
	\end{tabularx}
	\caption{Characteristics of the optimization problems in the firefighting literature.}
	\label{tab:review_characteristics}
\end{table}

The dynamic nature of the wildfire is addressed in different ways. \cite{Mendes_Alvelos} and \cite{Mendes_Alvelos_2} consider that the wildfire can move through the different nodes by moving between adjacent nodes. These frameworks extend the work of \cite{Alvelos2018}, in which the author proposes several mixed integer programming (MIP) models to optimize the location of extinguishing resources with different objective functions. In order to solve large data instances, \cite{Mendes_Alvelos} use an Iterated Local Search (ILS) metaheuristic, whereas \cite{Mendes_Alvelos_2} develop a Robust Tabu Search algorithm.

\cite{Belval} and \cite{Zhou} represent the terrain in a grid where each cell represents a parcel of land of homogeneous terrain. The wildfire moves between these cells taking into account the fire spread paths, which depend on weather scenarios. In \cite{Belval_2}, the authors add to this model fireline construction constraints for ground resources.

In \cite{Wiitala} the dynamism of the wildfire is interpreted by means of the amount of burned area that varies over time. In a similar way, \cite{Donovan}, \cite{Jorge_1} and \cite{Jorge_descomposicion} model the advance of the wildfire with an increase in the perimeter of the wildfire at each time interval. Another approach is that of \cite{Ntaimo}, which uses a wildfire simulator to monitor the progress of the wildfire.

The most similar representation of wildfire dynamism to the one chosen in this work is found in \cite{Skorin}, where the authors represent some water targets in the different zones of the wildfire that must be satisfied. These water targets vary over time, so that they adapt to the needs of the wildfire as it progresses.

\begin{table}[h]
	\centering
	\footnotesize
	\begin{tabularx}{\textwidth}{ccXc}
		Reference &  Model & Objectives & Methods  \\ 
		\hline	
		\rule{0pt}{3ex} \cite{Mendes_Alvelos} & MILP & Minimize the area burned and resources deployed & ILS, Gurobi \\
		\cite{Mendes_Alvelos_2} & MILP & Minimize the area burned and resources deployed & TS, Gurobi  \\
		\cite{Belval} & MILP & Minimize the area burned & CPLEX \\
		\cite{Wiitala} & NLMIP & Minimize costs & Dynamic programming \\
		\cite{Ntaimo} & MILP & Maximize wildfires contained and minimize costs & CPLEX \\
		\cite{Donovan} & MILP & Minimize costs & CPLEX\\
		\cite{Zhou} & MILP & Minimize costs \& people evacuated & Goal programming \\
		\cite{Jorge_1} & MILP & Minimize costs & CPLEX \\
		\cite{Jorge_descomposicion} & MILP & Minimize costs & Decomposition techniques \\
		\cite{shahidi} & MILP & Minimize extinguishing time & G, CPLEX\\
		\cite{Skorin} & MILP & Reach a water threshold and minimize output water & RG, SA, CPLEX \\
		Present work & MILP & Maximize water drops efficiency and minimize penalizations & SA, ILS, Gurobi \\
		\hline
		\multicolumn{4}{c}{ILS: Iterated Local Search, TS: Tabu Search, G: Greedy, RG: Random Greedy, SA: Simulated Annealing.} \\
	\end{tabularx}
	\caption{Characteristics of the optimization models and methods of resolution in firefighting literature.}
	\label{tab:references}
\end{table}

Table \ref{tab:references} presents models developed, their objectives and methods chosen to solve them in the mentioned frameworks. In general, when it is necessary to solve large data instances, the authors have to appeal to resolution methods, since they are not able to solve the problems with commercial solvers. 

The work developed in \cite{Jorge_1} has coincidences with what is proposed here, since they also perform an optimal planning of firefighting resources, modeling the work of the aircraft during the whole planning and taking into account the legislation. However, following the model proposed in \cite{Donovan}, the optimization criteria of the model are only economic. In order to solve big data instances, in \cite{Jorge_descomposicion} use decomposition techniques. In the present work other criteria are used, since currently in Spain economic criteria do not govern aircraft planning. Furthermore, we include the process of loading and dropping water of helicopters. 

The model presented in \cite{shahidi} is also similar to the one proposed in this paper. The authors formulate a Vehicle Routing Problem (VRP) with the aim of sending all the extinguishing resources to the points of the wildfire, minimizing the total extinguishing time. As in the model presented in this framework, they use a graph where the nodes are the wildfire points, the water loading points, and the starting points of the resources. Unlike present work, they consider water loading points to be infinite. Another important difference is that the demands assigned to each wildfire node are static, so the evolution of the wildfire is not taken into account. 

The model presented in \cite{Skorin} attempts to establish a planning of aerial extinguishing resources, with the primary objective of reaching a set water threshold in the different areas of the wildfire at each time period. Similarly to the present work, the authors take into account that aircraft work in groups and must comply with current aviation legislation. Another similarity is that the initial condition of the aircraft is also considered. However, there are some differences between both approaches. In fact, \cite{Skorin} assume certain simplifications. They assume that each wildfire zone has its associated water loading point and use a ratio of the number of drops that each aircraft can perform per time period. Moreover, their model assigns the maximum consecutive time as working time for any active aircraft. Similarly, once it starts resting, it will rest for the minimum mandatory time. To solve large data instances, a metaheuristic,  that combines a randomized greedy procedure with the simulated annealing technique, is considered.
 
Thus, to the best of our knowledge, there was no model in the literature covering all the issues to be addressed in the optimization problem studied in this framework. So, the aim of this framework is to obtain a decision tool to provide the wildfire coordinator with a feasible and reasonable planning of all the helicopter operations in any scenario. The model of this framework can be solved in small instances with a commercial solver. However, the model is too complex to find reasonable solutions with the commercial solver in affordable time for real-life situations. To tackle this question, we develop two algorithms based on metaheuristic techniques, that have shown very promising results in all scenarios. In both algorithms the initial solution is obtained by means of a greedy metaheuristic, where the main trajectories and the working times of the helicopters are determined. In the improving phase, both algorithms differ, since one of them makes use of the Simulated Annealing approach to improve the solution, while the other is based on the Iterated Local Search procedure.

The structure of this paper is as follows. In Chapter 2 all the essential elements of the problem to be tackled are described. Chapter 3 deals with the mathematical formulation of the model. In Chapter 4 the two algorithms based on metaheuristic techniques are explained. Finally, Chapter 5 is devoted to the analysis of the results and the drawing of conclusions.

\section{Problem Description}
The purpose of this paper is to provide a tool that aids decision-making for an aerial coordinator during a wildfire. 

When a wildfire occurs, it is almost always necessary to use aerial resources. Airplanes and helicopters cannot operate simultaneously in the same operational zone due to air safety reasons.
Helicopters and airplanes have distinct operating methods because, whereas the helicopters are capable of tighter maneuvers and hovering, airplanes rely on integrated tanks and require large water sources for refilling. For this reason, the  single use of helicopters is very common in firefighting tasks. 


Helicopters are equipped with Bambi buckets 
to collect water from standing sources, but loading from rivers requires specific conditions to avoid hazards. Helicopters typically follow an intuitive operating pattern, adhering to aviation regulations regarding flight hours and rest periods. The aerial coordinator plays a vital role in coordinating and planning all aircraft involved in firefighting operations. In order to plan the work of firefighting helicopters, the aerial
coordinator must make several decisions, all of which are interrelated: the helicopters that will operate in the firefighting, the working periods, the loading water points, the dropping areas and the rest bases.

Due to the complex nature of decision-making involved in helicopter operations, our focus in this work will be on providing a tool specifically designed for planning the utilization of these aircraft.

The main elements to consider in this decision problem are briefly described below.
\subsection{Wildfire}
\label{Wildfire}

\begin{figure}[!h]
	\begin{minipage}[c]{0.45\linewidth}
		\caption{Image obtained from the fire simulator FlamMap.}
		\centering
		\includegraphics[width=0.8\textwidth]{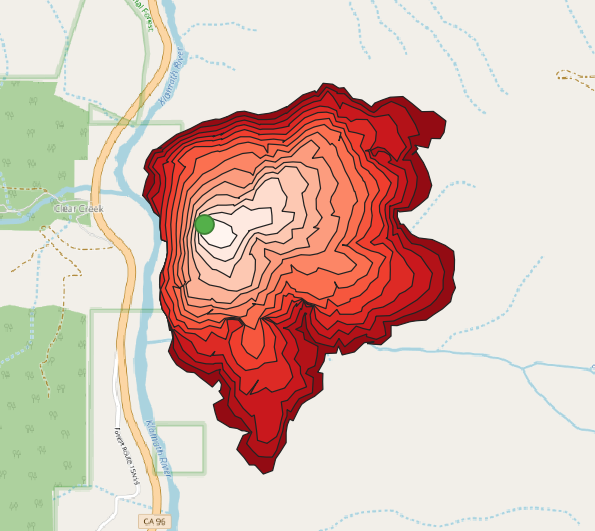}
		
		\label{fig:flammap_output}
	\end{minipage} \hfill
	\begin{minipage}[c]{0.45\linewidth} 
		\caption{Diagram used by wildfire coordinator representing the fire.} 
		\centering
		\includegraphics[width=0.65\textwidth]{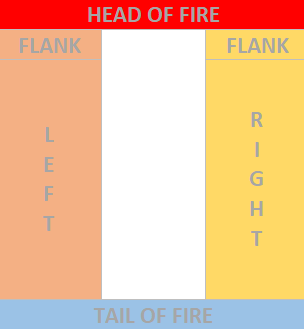}
		
		\label{fig:fire_scheme}
	\end{minipage}
\end{figure}

Once the aerial coordinator is on-site at the wildfire, understanding the anticipated progression of the wildfire is crucial. It can be determined using wildfire simulators or taking into account the coordinator's knowledge of terrain slope, vegetation, and expected wind patterns. Knowing the wildfire's evolution allows the coordinator to determine where aircraft should drop water, considering variations in temperature and wildfire spread across different zones. For visualization, Figure \ref{fig:flammap_output} illustrates an output of the FlamMap (\cite{FlamMap}) wildfire simulator, showing a wildfire's evolution from its point of ignition, where the green dot represents the starting point and the black lines depict the wildfire perimeter at intervals (e.g., every 30 minutes). The colors indicate the passage of time, where the darker shades indicate the expected perimeters at the end of the planning period, typically 7 hours after the wildfire's ignition.

The tools provide the wildfire coordinator with information on potential working areas, indicating where water drops would be most effective. They also help to define firefighter working areas, which need to meet certain safety criteria, such as easy accessibility and emergency exit routes. Aircraft work zones should be different from firefighters' areas for safety reasons, except in occasional instances when a firefighter brigade requests water drops from aircraft for reinforcement. Typically, firefighters work at the tail of a wildfire, while aircraft usually drop water on the head and flanks. The coordinator uses a diagram, as shown in Figure \ref{fig:fire_scheme}, to represent the wildfire and the aircraft action zones. In this diagram, the aerial coordinator marks the aircraft attack zones. These working areas may change over time, particularly for large wildfires that last several hours.

One of the most critical decisions for the aerial coordinator is determining the water drop zones for the aircraft, based on the expected wildfire evolution. Therefore, in the developed model, the area where each helicopter should drop water at each time interval is selected. To make this decision, the aerial coordinator must indicate the efficiency of water drops for each wildfire point, using the scheme shown in Figure \ref{fig:fire_scheme}. 
For each time interval, the coordinator assigns a value from 0 to 10 to each wildfire zone, with a higher value representing more efficient water drops, and 0 indicating no effect.

As the wildfire evolves over time, the efficiency values assigned to different areas will change accordingly. Initially, the aerial coordinator might prefer aircraft to drop water on the flanks of the wildfire but, after an hour, it might be more effective to target the head of the wildfire. 



\begin{figure}[H]
	\caption{Assigned efficiency to wildfire nodes over time, representing the progress of the wildfire.}
	\centering
	\includegraphics[width=0.85\textwidth]{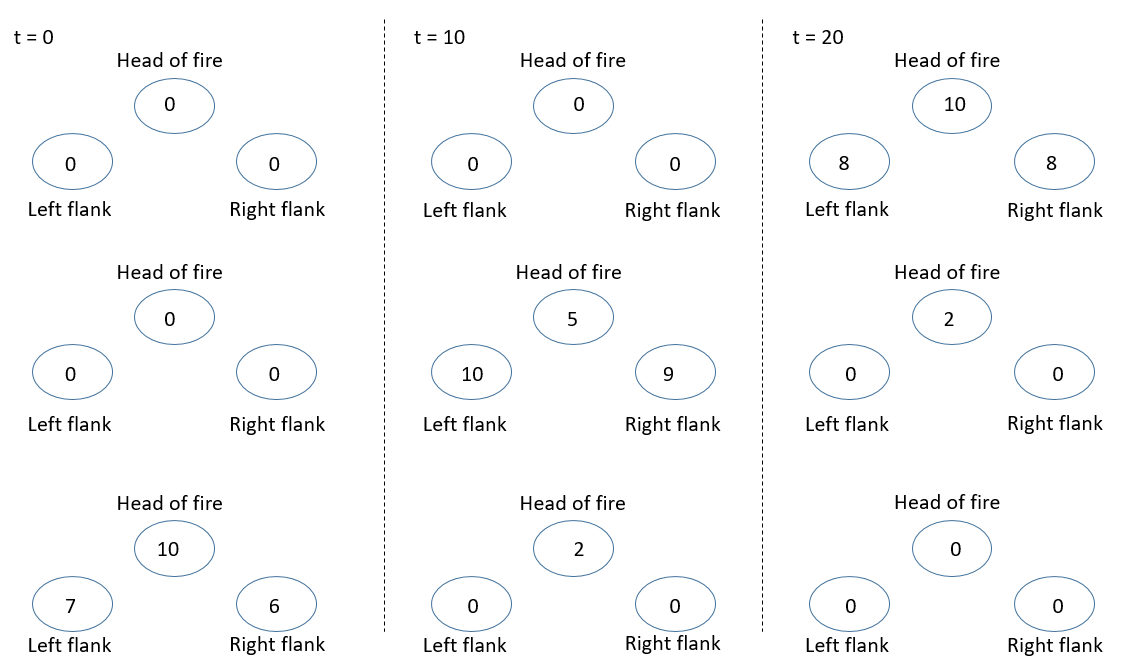}
	 
	\label{fig:ef_efficiency_expan}
\end{figure}

The evolution of the wildfire is represented by nodes which include expected future positions. For example, if the wildfire is predicted to move 2 km north, it can be represented by stages: the initial position, 1 km north after a certain time, and 2 km to the north at the end of the planning period. Figure \ref{fig:ef_efficiency_expan} shows the simulator output providing these positions and times. Consequently, 
the wildfire is represented by 9 nodes (three for each stage, representing the head and the two flanks of the wildfire). 
Initially, nodes representing future wildfire positions have an efficiency of 0, as there is no immediate need for water drops.

Figure \ref{fig:ef_efficiency_expan} illustrates this dynamic representation. Occasionally, the aerial coordinator might decide to preemptively drop water in areas where the wildfire has not yet reached, to cool the area and reduce potential flame intensity. These areas can be represented with new nodes, assigning the appropriate efficiency values to the relevant time intervals.

\subsection{Water load points}
Helicopters need to continuously load water while fighting wildfires. In Spain, all helicopters use a Bambi bucket system, 
which consists of a bucket suspended from the helicopter by a long cord, typically between 3.5 to 10 meters in length. To collect water, the helicopter hovers near standing water tanks where the Bambi bucket is submerged and loaded. These water tanks must be easily accessible, free from vegetation or obstacles that might hinder the helicopter's approach. Helicopters can also load water from rivers or the sea if necessary.

The aerial coordinator maintains a list of all possible water loading points, including safe river areas, reservoirs, lakes, and artificial loading points strategically built across Spain. These artificial points vary in capacity, and the coordinator knows the coordinates and capacities of all such tanks. In addition, private points, such as pools, can be used if needed, since their locations and capacities are also known by the coordinator.



\subsection{Rest bases}
Aircraft pilots need to take breaks from time to time, and use the opportunity to refuel their aircraft. This is done at the designated rest bases. Helicopters must arrive at the rest base empty of water, with the Bambi bucket folded. At the end of the planning all helicopters must be in a rest base.
In Spain, the “Circular Operativa 16-B” \cite{16Bravo} estates that helicopters should not fly for more than 2 consecutive hours, with a minimum break of 40 minutes.

The rest bases have a limited capacity to accommodate helicopters, depending on their surface area and refueling facilities. The aerial coordinator knows the coordinates and capacity of the rest areas in the region for which he/she is responsible.

In the region of Galicia, in 2023, according to the Galician Wildfire Prevention and Defense Plan (PLADIGA, \cite{Pladiga}), there were 20 helicopter rest bases distributed throughout the territory. 

\subsection{Helicopter operation}

When a helicopter is selected by the aerial coordinator to work on a wildfire, it flies to the wildfire area as soon as it receives the call. If the helicopter is resting after working on a previous wildfire, it must complete its rest before returning to work. If the helicopter is currently working on another wildfire, it will fly directly to the new wildfire if it has flight time enough; otherwise, it will rest at the nearest base before joining the new wildfire. Once at the wildfire, the helicopter joins the designated circuit, first flying to the water loading point and then to the indicated drop area.

Due to low visibility when flying over wildfires, helicopters follow established circuits called ``main trajectories'', flying one behind the other and maintaining a safe distance to avoid collisions. These circuits are closed loops between the water loading and dropping points. The aerial coordinator can establish multiple main trajectories, each with specific loading and dropping points, which remain constant until the wildfire's evolution changes.

The aerial coordinator decides the helicopters assigned to the main trajectories, and the loading and dropping points to be considered. At the beginning of a wildfire or in smaller fires, if the aerial coordinator is not present, the first helicopter decides the main trajectory, and the others follow it. Different classes of helicopters (heavy, medium, light) cannot share the same trajectory, though multiple trajectories of the same class can operate simultaneously. Once assigned, a helicopter remains in its group until it completes the assigned task.

\begin{figure}[H]
	\caption{Main trajectories on a wildfire occurred in Santiago de Compostela, Spain.} 
	\centering
	\includegraphics[width=0.5\textwidth]{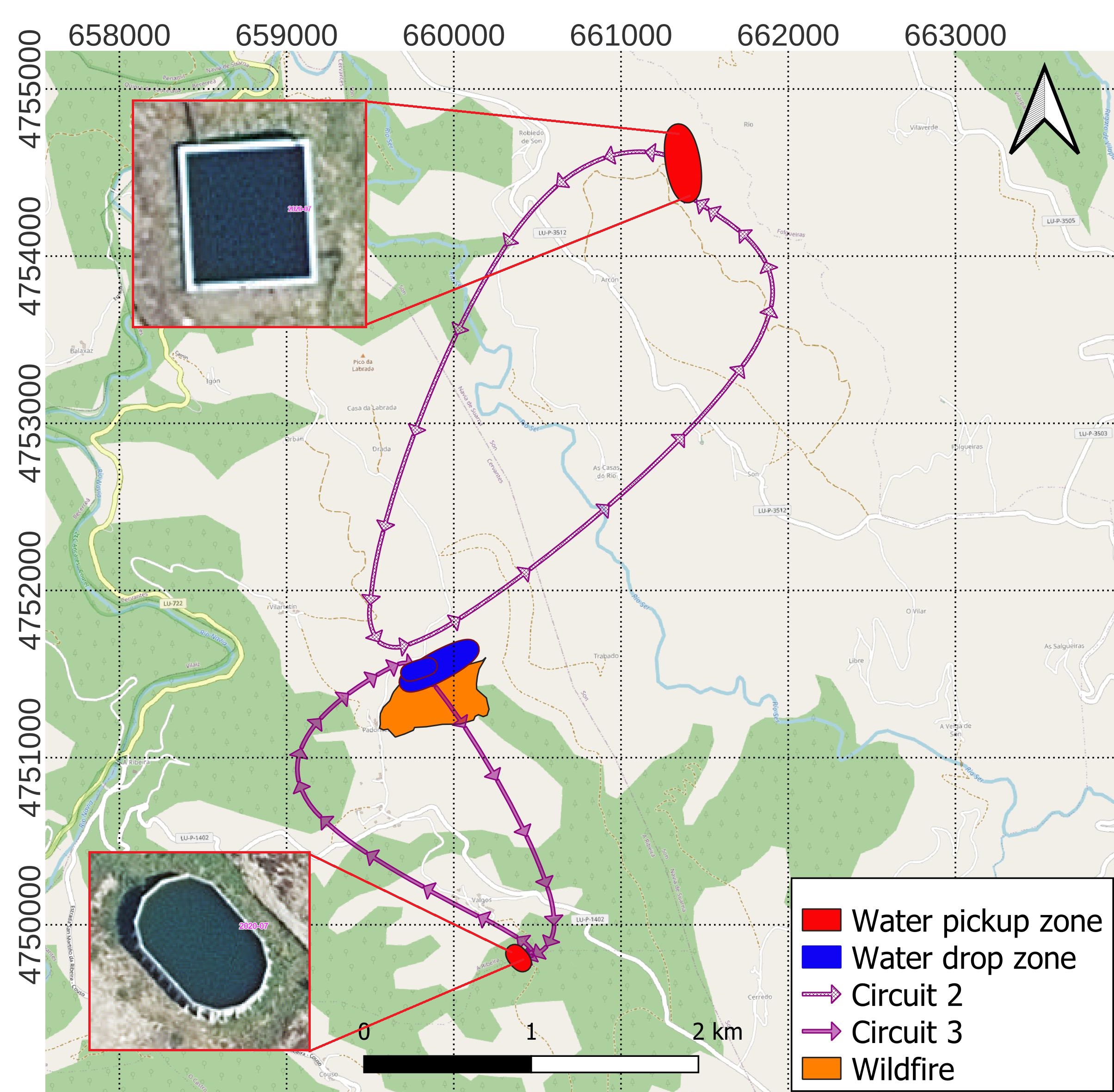}
	
	\label{fig:norias}
\end{figure}

Figure \ref{fig:norias} illustrates two main trajectories marked by arrows, which indicate the direction of the helicopters. The orange color represents the area of the wildfire. Red ellipses indicate the water loading areas, showing two artificial tanks used by the helicopters. Blue ellipses indicate the water drop zones. In this example, two main trajectories share the same drop zone but have different loading zones.

\section{Problem formulation}
The purpose of the model is to select adequately the helicopters operating on the wildfire, the loading water points, the drop zones and the rest bases. Thus, the aim of the model is to maximize the efficiency of helicopter water drops, taking into account current Spanish aviation legislation.

The helicopters operate in groups according to a specific main trajectory. The group to which each helicopter belongs is a parameter of the model, since it is a decision made by the aerial coordinator. Main trajectories are closed circuits that helicopters follow between the water loading point and the water dropping point. All helicopters belonging to the same main trajectory must follow the same circuit, that is, they must load and drop water at the same points. Each main trajectory is determined by its water loading and dropping zones, which remain unchanged until the wildfire changes. As the wildfire progresses, it is necessary for the helicopters to change their route to adapt to the new wildfire situation, so there may be a change in the water dropping zone or in the water loading zone.

If several helicopters of the same group are operating at the same time, they will fly behind each other, maintaining minimum safety distance. It means that two or more helicopters on the same main trajectory will never load or drop water at the same time.

\subsection{Mathematical model}
The firefighting scenario under study will be interpreted as a Mixed Integer Linear Programming (MILP) model. 

One of the main characteristics of the model is that the trajectories of the helicopters throughout the wildfire will be represented with the use of a graph $G = (N, E)$, where the set of nodes $N := N_A \cup N_C \cup N_I \cup N_B$ is composed of 4 types of nodes, indicating the initial positions of the aircraft ($N_A$), the water load points($N_C$), the wildfire points ($N_I$) and the rest bases ($N_B$). It is not necessary for the graph $G$ to be complete, since the set of arcs $E$ will only indicates the possible connections between any pair of nodes, if they exist. Moreover, the dynamic behavior of the wildfire will be implemented by means of the set of nodes, as are its positions, there are some associated parameters, such as the efficiency of the water drops in the wildfire points, that may vary over time.

Table \ref{tab:parameters} contains input parameters and decision variables of the model.

\begin{table}[h]
	\scriptsize
	\centering
	\begin{tabular}{ll}
		\hline
		\multicolumn{2}{l}{Sets} \\
		\hline
		$A$ & Set of helicopters. \\
		$W$ & Set of main trajectories.   \\
		$T$ & Set of time intervals.  \\
		$N = N_A \cup N_C \cup N_I \cup N_B$ & Set of nodes.\\
		$E$ & Set of edges. \\
		$N_A$ & Initial position of helicopters. \\
		$N_C$ & Water load points. \\
		$N_I$ & Wildfire points. \\
		$N_B$ & Rest bases.\\
		\hline
		\multicolumn{2}{l}{Input parameters} \\
		\hline
		$\lambda_{ija}$ & Time it takes to helicopter $a$ to fly up the edge $(i,j)$. \\
		$\alpha_{ia}$ & Time it takes to helicopter to load or drop water in the node $i$. \\
		${ef}_i^t$ & Water drop efficiency of node $i$ in time $t$. \\
		${CB}_i$ & Maximum number of helicopters which can rest at the same time in base $i$. \\
		${CA}_i$ & Initial water capacity (liters) of the water point $i$. \\
		${CL}_i$ & Number of helicopters which can load water at the same time in the point $i$. \\
		$n_{aw}$ & Takes value 1 if helicopter $a$ belongs to main trajectory $w$, 0 otherwise. \\
		$sp_a$ & Initial position of helicopter $a$, which must be a node of $N_A$. \\
		$wc_a$& Water capacity (liters) of helicopter $a$.\\
		$wli_a$& Initial water load status of helicopter $a$.\\
		$cfi_a$ & Consecutive flight time of helicopter $a$ when the model is executed.\\
		$tf_a$ & Total flight time of helicopter $a$ during the day of model execution.\\
		$mcf_a$& Maximum consecutive flying time of helicopter $a$.\\
		$mtf_a$& Maximum total flying time of helicopter $a$.\\
		$ri_a$& Consecutive rest time of helicopter $a$ when the model is executed.\\
		$mr_a$& Minimum consecutive rest time to be performed by aircraft $a\in A$.\\
		$ev^t$& Takes value 1 if in the time interval $t\in T$ there is a wildfire evolution.\\
		\hline
		\multicolumn{2}{l}{Decision variables} \\
		\hline
		$y_{ia}^t \in \{0,1\}$ & Takes value 1 if helicopter $a$ is at node $i$ in $t$, 0 otherwise.\\
		$x_{ija}^t \in \{0,1\}$& Takes value 1 if helicopter $a$ flies over edge $(i,j)$ in $t$, 0 otherwise.\\
		$e_{ia}^t \in \{0,1\}$& Takes value 1 if helicopter $a$ ends a break at rest base $i$ in $t$, 0 otherwise.\\ 
		$ec_{ia}^t \in \{0,1\}$& Takes value 1 if helicopter $a\in A$ finishes loading water at node $i \in N_C$ in $t\in T$, \\ & 0 otherwise.\\  
		 $ed_{ia}^t \in \{0,1\}$& Takes value 1 if helicopter $a$ finishes dropping water at node $i$ in $t$, 0 otherwise.\\  
		 $r_{iw}^t \in \{0,1\}$&  Takes value 1 if main trajectory $w$ is associated with node $i$ in $t$, 0 otherwise.\\ 
		 $cw_w^t \in \{0,1\}$& Takes value 1 if main trajectory $w$ can change its associated water or dropping zone in\\ &  $t$,  0 otherwise.\\ 
		 $aux_w^t \in \{0,1\}$& Binary variable that shows if main trajectory $w$ makes a change in $t$.\\ 
		 $h_1^t \in \mathbb{R}^{+}$& Slack variable to prevent infeasibilities  in $t$.\\ 
		 $faux_a^t \in \mathbb{R}^{+}$& Counts the number of intervals remaining for helicopter $a$ to work the maximum  \\ & number of consecutive  intervals before the break in $t$.\\ 
		 \hline
	\end{tabular}
	\caption{Input parameters and decision variables of the model.}
	\label{tab:parameters}
\end{table}

Furthermore, there are a number of auxiliary variables, which are detailed below.

\begin{equation}
	cnta_{i}^t := \sum_{a\in A} \sum_{t'=1}^t y_{ia}^t cap_a \;\;\; \forall i\in N_C, t\in T 
\end{equation}
\begin{equation}
	ca _{i}^t := CA_i - cnta_i^{t-1} \;\;\; \forall i\in N_C, t>1  \end{equation}
\begin{equation}ca _{i}^1 := CA_i  \;\;\; \forall i\in N_C  
\end{equation}
\begin{equation} z_a^t := z_a^{t-1} + \sum_{i\in N_c}ec_{ia}^t - \sum_{j\in N_I}ed_{ia}^t \;\;\; \forall a\in A, \; t>1
\end{equation}
\begin{equation} z_a^1 := wli_a \;\;\; \forall a\in A
\end{equation} 
\begin{equation}cntv_a^t := cfi_a + \sum_{(i,j)\in E}\sum_{l=1}^{t}x_{ija}^l + \sum_{i\in N_C \cup N_I} \sum_{l=1}^{t} y_{ia}^l  + faux_a^t - \sum_{l=1}^{t}\sum_{i\in N_B}mcf_a e_{ia}^l \;\;\; \forall t \in T, a\in A
\end{equation}

(1) is a water counter of each water load point. (2) and (3) are the updated capacity of each water load point. (4) is an indicator of aircraft water load status, takes value 1 if helicopter $a$ is loaded during the time interval $t$, 0 otherwise. (5) is a time counter of consecutive flight intervals of each helicopter.

\subsubsection{Objective function}
The objective function consists of maximizing the efficiency associated with the water drops and minimizing a penalty associated with the flights of the helicopters, the changes in the main trajectories established and the number of time intervals in which no helicopters are working to extinguish the wildfire. 

\begin{equation*}
	\text{Max}\; (\sum_{t\in T} \sum_{i\in N_I} \sum_{a\in A} y_{ia}^t {ef}_i^t wc_a)/ub_1  - \mu_2(\sum_{t\in T} \sum_{(i,j)\in E} \sum_{a\in A} x_{ija}^t)/ub_2 - $$ $$ - \mu_3(\sum_{t\in T} \sum_{i \in N_C \cup N_I} \sum_{a\in A}  y_{ia}^t)/ub_3  -  \mu_4(\sum_{t\in T}\sum_{w \in W}aux_w^t)/ub_4 - \mu_5\sum_{t\in T}h_1^t - $$ $$ - \mu_6\sum_{t\in T}\sum_{a\in A}faux_a^t,
\end{equation*}

where $\mu_i$, $i=2,\cdots,6$ parameters are the normalized weights assigned to each term. The fundamental objective of the problem is to maximize the efficiency of water drops. The weights assigned to each term of the objective function are chosen to take into account the different objectives of the model in lexicographic order ($\mu_2 > \mu_3 > \mu_4 > \mu_5 > \mu_6$).

\setcounter{equation}{0}
\subsubsection{Model Constraints}
\begin{equation}
	\sum_{a \in A} y_{ia}^t \leq {CB}_i \;\;\; \forall t\in T, i \in N_B
\end{equation}
\begin{equation} \sum_{a \in A} y_{ia}^t \leq CL_i \;\;\; \forall t\in T, i \in N_C
\end{equation}
\begin{equation} \sum_{a \in A} y_{ia}^t wc_a \leq ca_i^t \;\;\; \forall t>1, i \in N_C
\end{equation}
\begin{equation} \sum_{a \in A} y_{ia}^1 wc_a \leq CA_i \;\;\; \forall  i \in N_C
\end{equation}
\begin{equation} y_{ia}^t  \leq z_a^t \;\;\; \forall t\in T, i \in N_I, a\in A
\end{equation}
\begin{equation} y_{ia}^t  \leq 1 - z_a^t \;\;\; \forall t\in T, i \in N_C, a\in A
\end{equation}
\begin{equation}  y_{ia}^t - y_{ia}^{t-1} \leq y_{ia}^l \;\;\; \forall i\in N_C \cup N_I, a\in A,t>1, l\in \{t+1,...,min\{|T|,t+\alpha_{ia}-1\}\} / \alpha_{ia}>1
\end{equation}
\begin{equation} y_{ia}^t - y_{ia}^{t-1} \leq ed_{ia}^{min\{|T|,t +\alpha_{ia}\}} \;\;\; \forall i\in N_I, a\in A, t>1
\end{equation}
\begin{equation} y_{ia}^t - y_{ia}^{t-1} \leq ec_{ia}^{min\{|T|,t +\alpha_{ia}\}} \;\;\; \forall i\in N_C, a\in A, t>1
\end{equation}	
\begin{equation} y_{ia}^{t-1} - y_{ia}^t  \leq e_{ia}^{t-1} \;\;\; \forall i\in N_B, a\in A, t>1
\end{equation}
\begin{equation}e_{ia}^t \leq y_{ia}^{t} \;\;\; \forall b\in N_B, a\in A
\end{equation}
\begin{equation} y_{ia}^t \leq 1 - e_{ia}^{t-1} \;\;\; \forall i\in N_B, a\in A, t>1
\end{equation}
\begin{equation} y_{ia}^{t-1} - y_{ia}^t  \leq ec_{ia}^t \;\;\; \forall i\in N_C, a\in A, t>1
\end{equation}
\begin{equation}ec_{ia}^t \leq y_{ia}^{t-1} \;\;\; \forall i\in N_C, a\in A, t>1
\end{equation}
\begin{equation} y_{ia}^t \leq 1 - ec_{ia}^t \;\;\; \forall i\in N_C, a\in A, t>1
\end{equation}
\begin{equation} y_{ia}^{t-1} - y_{ia}^t  \leq ed_{ia}^t \;\;\; \forall i\in N_I, a\in A, t>1
\end{equation}
\begin{equation}ed_{ia}^t \leq y_{ia}^{t-1} \;\;\; \forall i\in N_I, a\in A, t>1
\end{equation}
\begin{equation} y_{ia}^t \leq 1 - ed_{ia}^t \;\;\; \forall i\in N_I, a\in A, t>1
\end{equation}	
\begin{equation}tf_a + \sum_{(i,j) \in E} \sum_{t \in T} x_{ija}^t  + \sum_{t \in T} \sum_{i \in N_C \cup N_I} y_{ia}^t  \leq mtf_a \;\;\; \forall a \in A
\end{equation} 
\begin{equation}cntv_a^t \leq mcf_a \;\;\; \forall t\in T, a\in A 
\end{equation} 
\begin{equation}cntv_a^t \geq 0 \;\;\; \forall t\in T, a\in A 
\end{equation} 
\begin{equation} \sum_{l = t}^{min\{|T|,t+mr_a-1\}} y_{ia}^l \geq mr_a - mr_a(1-y_{ia}^t+y_{ia}^{t-1}-e_{ia}^{t-1}c) \;\;\; \forall i\in N_B, a\in A, t>1 
\end{equation} 
\begin{equation} \sum_{l = 1}^{min\{|T|,mr_a-ri_a\}} y_{ia}^l \geq mr_a - ri_a \;\;\; \forall i\in N_I, a\in A, mr_a>ri_a>0 
\end{equation}
\begin{equation}  \sum_{i\in N_B}y_{ia}^{|T|} \geq 1 \;\;\; \forall a\in A 
\end{equation}
\begin{equation}y_{ia}^1 = 1 \;\;\; \forall a\in A, \; i\in sp_a
\end{equation}	
\begin{equation} y_{ja}^{min\{|T|,t+\lambda_{ija}\}} \geq x_{ija}^t - x_{ija}^{t-1} \;\;\; \forall (i,j)\in E, a\in A, t>1 
\end{equation}	.
\begin{equation} x_{ija}^t - x_{ija}^{t-1} \leq y_{ia}^{t-1} \;\;\; \forall (i,j)\in E, a\in A, t>1 
\end{equation}
\begin{equation}  y_{ja}^t \leq \sum_{i \in N / (i,j)\in E}x_{ija}^{t-1} + y_{ja}^{t-1} \;\;\; \forall j\in N\setminus N_A, a\in A, t>1   
\end{equation}
\begin{equation}  y_{ia}^t - y_{ia}^{t+1} \leq \sum_{j \in N / (i,j)\in E}x_{ija}^{t+1} \;\;\; \forall i\in N, a\in A, t<|T|  
\end{equation}
\begin{equation} x_{ija}^t - x_{ija}^{t-1} \leq x_{ija}^l \;\;\; \forall (i,j)\in E, a\in A,t>1, l\in \{t,...,min\{|T|,t+\lambda_{ija}-1\}\}
\end{equation}
\begin{equation}\sum_{i\in N/(i,j)\in E}x_{ija}^t \leq \sum_{i\in N/(j,i)\in E}\sum_{l > t}x_{jia}^l \;\;\; \forall j\in N/(N_B \cup N_A), \; t,l\in T, a\in A
\end{equation}
\begin{equation}\sum_{i\in N/(i,j)\in E}x_{ija}^t \leq \sum_{i\in N/(j,i)\in E}\sum_{l > t}x_{jia}^l + y_{ja}^l \;\;\; \forall j\in N_B,\; t,l\in T,  a\in A
\end{equation}
\begin{equation}\sum_{i\in N} y_{ia}^t + \sum_{(i,j) \in E}x_{ija}^t \leq 1 \;\;\; \forall a\in A, \; t\in T
\end{equation}
\begin{equation}\sum_{i\in N_C \cup N_I} r_{iw}^t \leq 1 \;\;\; \forall w\in W, \; t\in T
\end{equation}
\begin{equation}r_{iw}^{t-1} \leq r_{iw}^{t} + cw_w^{t} \;\;\; \forall i\in N_C\cup N_I, \; w\in W,\; t>1
\end{equation}
\begin{equation}aux_w^t \geq r_{iw}^t - r_{iw}^{t-1} \;\;\; \forall i\in N_C\cup N_I, \;w\in W, \; t>1
\end{equation}	
\begin{equation}cw_w^t \leq ev^t + cw_w^{t-1}- aux_w^{t-1} \;\;\; \forall w\in W, \; t>1
\end{equation}
\begin{equation}cw_w^t \geq cw_w^{t-1} -  aux_w^{t-1} \;\;\; \forall w\in W, \; t>1
\end{equation}
\begin{equation}cw_w^t \geq ev^t \;\;\; \forall w\in W, \; t>1
\end{equation}
\begin{equation} \sum_{a\in A} y_{ia}^t n_{aw} \leq 1 \;\;\; \forall i\in N_C \cup N_I, w\in W, t\in T
\end{equation}
\begin{equation}y_{ia}^t n_{aw} \leq r_{iw}^t \;\;\; \forall i\in N_C\cup N_I, \;w\in W, \; a\in A, \; t\in T
\end{equation}
\begin{equation}\sum_{(i,j)\in E}\sum_{k\in N_C \cup N_I}\sum_{a\in A}x_{ija}^t+y_{ka}^t+h_1^t \geq 1 \;\;\; \forall t\in T,\; |T| >t>1
\end{equation}
\begin{equation}
	faux_a^t >= faux_a^{t-1}  \;\;\; \forall a\in A, \; t\in T,\; t > 1
\end{equation}
\begin{equation} e_{ia}^1 := 0 \;\;\; \forall i\in N_B, a\in A
\end{equation}
\begin{equation} ec_{ia}^1 := 0 \;\;\; \forall i\in N_C, a\in A
\end{equation}
\begin{equation} ed_{ia}^1 := 0 \;\;\; \forall i\in N_I, a\in A
\end{equation}
\begin{equation}aux_w^1 := 0 \;\;\; \forall w\in W
\end{equation}
\begin{equation}cw_w^1 := 0 \;\;\; \forall w\in W
\end{equation}

Constraint (1) prevents the capacities of the resting bases from being exceeded.
Constraints (2)-(4) avoid exceeding the capacities of the water loading points, both in terms of the number of helicopters loading water at each interval of time and the amount of water extracted from each point. In Constraints (5)-(18) some relationships between variables are reflected. For its part, Constraints (5)-(6) model the relationship between the state of the helicopter's water load and the action of loading or dropping it. The intervals of time when the aircraft is loading/dropping water or leaving the corresponding node are indicated by (7)-(9). Constraints (10)-(12) are related to the intervals of time at which helicopters start or end its rests. In a similar way (13)-(18) are associated with the intervals of time at which helicopters load or drop water. Constraints (19)-(23) ensure compliance with flight regulations. Constraint (19) guarantees that the maximum daily flight time is not exceeded. Constraints (20)-(21) prevent the maximum consecutive working time from being exceeded, whereas (22)-(23) regulate mandatory rest breaks. Constraints (24)-(25) indicate the helicopters positions at the beginning and end of planning. The helicopters movements in the graph are controlled by (26)-(33). The constraints (34)-(41) are related to the main trajectories. The right behavior of helicopters in the same trajectory is controlled by (40)-(41), since they cannot load or drop water at the same period and must use the load or drop points of the trajectory they belong to. Constraint (42) models that whenever possible at least one helicopter should be working on the wildfire. Constraint (43) increases the value of the helicopters flight interval counter. Finally, Constraints (44)-(48) are devoted to set the initial values of some variables.

To find the optimal solutions in this model for small instances, the commercial solver Gurobi has been applied. However, as it will be shown in Section 5, it is very difficult to find competitive solutions (or even feasible solutions) with the solver in realistic instances. Thus, Section 4 deals with two new proposals based on metaheuristic techniques.

%
%
%

\section{Metaheuristic algorithms}
In view of the emergency situation tackled in the model proposed in the previous section, a sufficiently good feasible solution should be found in a relatively short time, preferably, a maximum of 15 minutes. Since commercial solver Gurobi has not been capable of finding optimal solutions in satisfactory computational times (and in some cases not even feasible solutions), it was necessary to design \textit{ad hoc} algorithms adapted to the characteristics of the problem.

To this aim, two algorithms (from now on G+SA and G+ILS) were developed. Both algorithms are based on metaheuristic techniques tailored to the characteristics of the problem under study.

The first stage of the algorithms, which corresponds to the construction of a feasible initial planning, is similar in both algorithms. To obtain this feasible planning, a Greedy procedure was proposed. The second stage, devoted to the improvement of the solution, differs in both algorithms. One of the algorithms (G+SA) applies the Simulated Annealing technique (\cite{kirkpatrick}), while the other (G+ILS) makes use of the Iterated Local Search method (\cite{handbook_met2}).

In the following, the steps to obtain the initial solution are described in detail, as well as the solution improvement phase in both algorithms.

\subsection{Initial solution} \label{sec:initial_solution}
The Greedy procedure to obtain the initial solution consists of two parts: in the first stage, the main trajectories are defined by identifying their load and water drop points\footnote{Notice that the number of main trajectories and the helicopters belonging to each main trajectory are inputs of the problem, since all of them are previously defined by the aerial coordinator.}. In the second stage, the working time for all the aircraft of each main trajectory is determined.

\begin{figure}[h]
	\caption{Scheme of the process of defining main trajectories for the initial solution.}
\tikzset{every picture/.style={line width=0.75pt}} 
\begin{tikzpicture}[x=0.75pt,y=0.75pt,yscale=-1,xscale=1]
	
	\draw   (36,59.6) -- (205,59.6) -- (205,105.6) -- (36,105.6) -- cycle ;
	\draw   (113,25.6) .. controls (113,18.42) and (123.75,12.6) .. (137,12.6) .. controls (150.25,12.6) and (161,18.42) .. (161,25.6) .. controls (161,32.78) and (150.25,38.6) .. (137,38.6) .. controls (123.75,38.6) and (113,32.78) .. (113,25.6) -- cycle ;
	\draw  [dash pattern={on 0.84pt off 2.51pt}]  (101,25.6) .. controls (72.14,22.61) and (-34.92,55.27) .. (29.02,89.48) ;
	\draw [shift={(30,90)}, rotate = 207.53] [color={rgb, 255:red, 0; green, 0; blue, 0 }  ][line width=0.75]    (10.93,-3.29) .. controls (6.95,-1.4) and (3.31,-0.3) .. (0,0) .. controls (3.31,0.3) and (6.95,1.4) .. (10.93,3.29)   ;
	\draw   (303,58.6) -- (381,58.6) -- (381,100.6) -- (303,100.6) -- cycle ;
	\draw    (211,82.6) -- (291,82.6) ;
	\draw [shift={(293,82.6)}, rotate = 180] [color={rgb, 255:red, 0; green, 0; blue, 0 }  ][line width=0.75]    (10.93,-3.29) .. controls (6.95,-1.4) and (3.31,-0.3) .. (0,0) .. controls (3.31,0.3) and (6.95,1.4) .. (10.93,3.29)   ;
	\draw   (257,174.6) -- (326,174.6) -- (326,212.6) -- (257,212.6) -- cycle ;
	\draw    (338,109.6) -- (304.97,168.85) ;
	\draw [shift={(304,170.6)}, rotate = 299.13] [color={rgb, 255:red, 0; green, 0; blue, 0 }  ][line width=0.75]    (10.93,-3.29) .. controls (6.95,-1.4) and (3.31,-0.3) .. (0,0) .. controls (3.31,0.3) and (6.95,1.4) .. (10.93,3.29)   ;
	\draw   (160,179.6) .. controls (160,176.29) and (162.69,173.6) .. (166,173.6) -- (186,173.6) .. controls (189.31,173.6) and (192,176.29) .. (192,179.6) -- (192,197.6) .. controls (192,200.91) and (189.31,203.6) .. (186,203.6) -- (166,203.6) .. controls (162.69,203.6) and (160,200.91) .. (160,197.6) -- cycle ;
	\draw   (373,182.8) .. controls (373,179.93) and (375.33,177.6) .. (378.2,177.6) -- (396.8,177.6) .. controls (399.67,177.6) and (402,179.93) .. (402,182.8) -- (402,198.4) .. controls (402,201.27) and (399.67,203.6) .. (396.8,203.6) -- (378.2,203.6) .. controls (375.33,203.6) and (373,201.27) .. (373,198.4) -- cycle ;
	\draw  [dash pattern={on 4.5pt off 4.5pt}]  (246,189.6) -- (207,189.6) ;
	\draw [shift={(205,189.6)}, rotate = 360] [color={rgb, 255:red, 0; green, 0; blue, 0 }  ][line width=0.75]    (10.93,-3.29) .. controls (6.95,-1.4) and (3.31,-0.3) .. (0,0) .. controls (3.31,0.3) and (6.95,1.4) .. (10.93,3.29)   ;
	\draw    (161,170.6) -- (136.8,115.43) ;
	\draw [shift={(136,113.6)}, rotate = 66.32] [color={rgb, 255:red, 0; green, 0; blue, 0 }  ][line width=0.75]    (10.93,-3.29) .. controls (6.95,-1.4) and (3.31,-0.3) .. (0,0) .. controls (3.31,0.3) and (6.95,1.4) .. (10.93,3.29)   ;
	\draw  [dash pattern={on 4.5pt off 4.5pt}] (364,131.1) .. controls (364,123.09) and (379.67,116.6) .. (399,116.6) .. controls (418.33,116.6) and (434,123.09) .. (434,131.1) .. controls (434,139.11) and (418.33,145.6) .. (399,145.6) .. controls (379.67,145.6) and (364,139.11) .. (364,131.1) -- cycle ;
	\draw  [dash pattern={on 0.84pt off 2.51pt}]  (358,130.6) -- (334,131.52) ;
	\draw [shift={(332,131.6)}, rotate = 357.8] [color={rgb, 255:red, 0; green, 0; blue, 0 }  ][line width=0.75]    (10.93,-3.29) .. controls (6.95,-1.4) and (3.31,-0.3) .. (0,0) .. controls (3.31,0.3) and (6.95,1.4) .. (10.93,3.29)   ;
	\draw  [dash pattern={on 4.5pt off 4.5pt}]  (337,189.6) -- (365,190.53) ;
	\draw [shift={(367,190.6)}, rotate = 181.91] [color={rgb, 255:red, 0; green, 0; blue, 0 }  ][line width=0.75]    (10.93,-3.29) .. controls (6.95,-1.4) and (3.31,-0.3) .. (0,0) .. controls (3.31,0.3) and (6.95,1.4) .. (10.93,3.29)   ;
	\draw    (392,170.6) -- (394.65,155.57) ;
	\draw [shift={(395,153.6)}, rotate = 100.01] [color={rgb, 255:red, 0; green, 0; blue, 0 }  ][line width=0.75]    (10.93,-3.29) .. controls (6.95,-1.4) and (3.31,-0.3) .. (0,0) .. controls (3.31,0.3) and (6.95,1.4) .. (10.93,3.29)   ;
	
	\draw (38,62.6) node [anchor=north west][inner sep=0.75pt]  [font=\footnotesize] [align=left] {wildfire point with the highest\\ water drop efficiency};
	\draw (122,19.4) node [anchor=north west][inner sep=0.75pt]  [font=\footnotesize]  {$t=0$};
	\draw (304,67) node [anchor=north west][inner sep=0.75pt]  [font=\footnotesize] [align=left] {closest water\\ load point};
	\draw (258,177.6) node [anchor=north west][inner sep=0.75pt]  [font=\footnotesize] [align=left] {wildfire\\ evolution?};
	\draw (162,182.6) node [anchor=north west][inner sep=0.75pt]  [font=\footnotesize] [align=left] {yes};
	\draw (379,182.8) node [anchor=north west][inner sep=0.75pt]  [font=\footnotesize] [align=left] {no};
	\draw (371,125) node [anchor=north west][inner sep=0.75pt]  [font=\footnotesize] [align=left] {increase $\displaystyle t$};	
\end{tikzpicture}	

\label{fig:diagramanorias}
\end{figure}
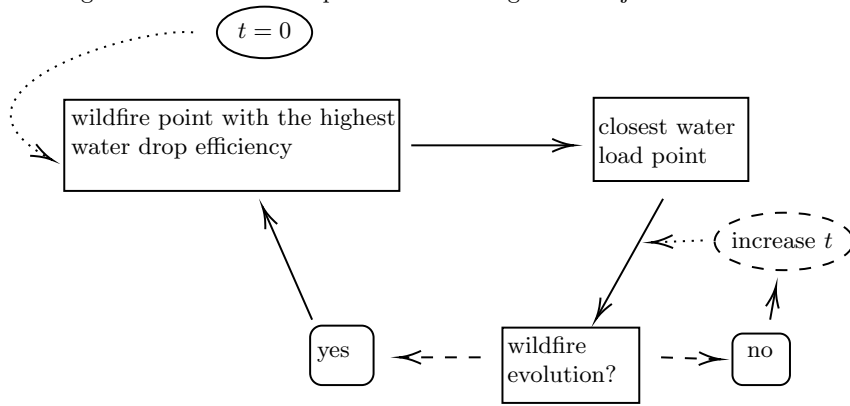

Figure \ref{fig:diagramanorias} shows an scheme of the definition of main trajectories. First, the wildfire point with the highest associated water drop efficiency is selected. Next, the closest water loading point is selected. In case there are more than one point with the same values, it is selected randomly. Both points are assigned to the main trajectory until there is a wildfire evolution, since all that time the efficiency values assigned to the wildfire might change.

In this way, each main trajectory has been assigned its water load and drop points for each time interval. Each helicopter belonging to a main trajectory must load or drop water at these points in the corresponding time interval.

For all main trajectories in the schedule, this process would be repeated. Then, for each main trajectory, the helicopters that must work at each time interval are selected.

\begin{figure}[h]
		\caption{Scheme of the helicopter work assignment process.}
	
	\tikzset{every picture/.style={line width=0.75pt}} 
	
	\begin{tikzpicture}[x=0.75pt,y=0.75pt,yscale=-1,xscale=1]
		
		\draw   (36,60.6) -- (150,60.6) -- (150,114.6) -- (36,114.6) -- cycle ;
		\draw   (63,29.1) .. controls (63,21.64) and (75.76,15.6) .. (91.5,15.6) .. controls (107.24,15.6) and (120,21.64) .. (120,29.1) .. controls (120,36.56) and (107.24,42.6) .. (91.5,42.6) .. controls (75.76,42.6) and (63,36.56) .. (63,29.1) -- cycle ;
		\draw  [dash pattern={on 0.84pt off 2.51pt}]  (58,20.6) .. controls (32.52,11.78) and (2.24,63.47) .. (25.51,82.48) ;
		\draw [shift={(27,83.6)}, rotate = 214.7] [color={rgb, 255:red, 0; green, 0; blue, 0 }  ][line width=0.75]    (10.93,-3.29) .. controls (6.95,-1.4) and (3.31,-0.3) .. (0,0) .. controls (3.31,0.3) and (6.95,1.4) .. (10.93,3.29)   ;
		\draw   (234,73.6) -- (302,73.6) -- (302,119.6) -- (234,119.6) -- cycle ;
		\draw   (401,98.1) .. controls (401,79.05) and (435.92,63.6) .. (479,63.6) .. controls (522.08,63.6) and (557,79.05) .. (557,98.1) .. controls (557,117.15) and (522.08,132.6) .. (479,132.6) .. controls (435.92,132.6) and (401,117.15) .. (401,98.1) -- cycle ;
		\draw   (308.4,24.8) .. controls (308.4,14.42) and (316.82,6) .. (327.2,6) .. controls (337.58,6) and (346,14.42) .. (346,24.8) .. controls (346,35.18) and (337.58,43.6) .. (327.2,43.6) .. controls (316.82,43.6) and (308.4,35.18) .. (308.4,24.8) -- cycle ;
		\draw   (376.4,159.8) .. controls (376.4,152.18) and (382.58,146) .. (390.2,146) .. controls (397.82,146) and (404,152.18) .. (404,159.8) .. controls (404,167.42) and (397.82,173.6) .. (390.2,173.6) .. controls (382.58,173.6) and (376.4,167.42) .. (376.4,159.8) -- cycle ;
		\draw  [dash pattern={on 4.5pt off 4.5pt}]  (438,63.6) -- (354.87,32.3) ;
		\draw [shift={(353,31.6)}, rotate = 20.63] [color={rgb, 255:red, 0; green, 0; blue, 0 }  ][line width=0.75]    (10.93,-3.29) .. controls (6.95,-1.4) and (3.31,-0.3) .. (0,0) .. controls (3.31,0.3) and (6.95,1.4) .. (10.93,3.29)   ;
		\draw  [dash pattern={on 4.5pt off 4.5pt}]  (480,127.6) -- (412.85,154.85) ;
		\draw [shift={(411,155.6)}, rotate = 337.91] [color={rgb, 255:red, 0; green, 0; blue, 0 }  ][line width=0.75]    (10.93,-3.29) .. controls (6.95,-1.4) and (3.31,-0.3) .. (0,0) .. controls (3.31,0.3) and (6.95,1.4) .. (10.93,3.29)   ;
		\draw    (298,30.6) -- (159.95,62.15) ;
		\draw [shift={(158,62.6)}, rotate = 347.12] [color={rgb, 255:red, 0; green, 0; blue, 0 }  ][line width=0.75]    (10.93,-3.29) .. controls (6.95,-1.4) and (3.31,-0.3) .. (0,0) .. controls (3.31,0.3) and (6.95,1.4) .. (10.93,3.29)   ;
		\draw    (159,94.6) -- (227,95.57) ;
		\draw [shift={(229,95.6)}, rotate = 180.82] [color={rgb, 255:red, 0; green, 0; blue, 0 }  ][line width=0.75]    (10.93,-3.29) .. controls (6.95,-1.4) and (3.31,-0.3) .. (0,0) .. controls (3.31,0.3) and (6.95,1.4) .. (10.93,3.29)   ;
		\draw    (310,98.6) -- (394,97.62) ;
		\draw [shift={(396,97.6)}, rotate = 179.33] [color={rgb, 255:red, 0; green, 0; blue, 0 }  ][line width=0.75]    (10.93,-3.29) .. controls (6.95,-1.4) and (3.31,-0.3) .. (0,0) .. controls (3.31,0.3) and (6.95,1.4) .. (10.93,3.29)   ;
		\draw   (336,200) -- (443,200) -- (443,235.6) -- (336,235.6) -- cycle ;
		\draw    (392.2,179.6) -- (392.91,195.6) ;
		\draw [shift={(393,197.6)}, rotate = 267.46] [color={rgb, 255:red, 0; green, 0; blue, 0 }  ][line width=0.75]    (10.93,-3.29) .. controls (6.95,-1.4) and (3.31,-0.3) .. (0,0) .. controls (3.31,0.3) and (6.95,1.4) .. (10.93,3.29)   ;
		\draw   (322,314.1) .. controls (322,295.05) and (356.03,279.6) .. (398,279.6) .. controls (439.97,279.6) and (474,295.05) .. (474,314.1) .. controls (474,333.15) and (439.97,348.6) .. (398,348.6) .. controls (356.03,348.6) and (322,333.15) .. (322,314.1) -- cycle ;
		\draw    (394,242.6) -- (394,270.6) ;
		\draw [shift={(394,272.6)}, rotate = 270] [color={rgb, 255:red, 0; green, 0; blue, 0 }  ][line width=0.75]    (10.93,-3.29) .. controls (6.95,-1.4) and (3.31,-0.3) .. (0,0) .. controls (3.31,0.3) and (6.95,1.4) .. (10.93,3.29)   ;
		\draw  [dash pattern={on 4.5pt off 4.5pt}]  (479,306.6) -- (509.15,294.35) ;
		\draw [shift={(511,293.6)}, rotate = 157.89] [color={rgb, 255:red, 0; green, 0; blue, 0 }  ][line width=0.75]    (10.93,-3.29) .. controls (6.95,-1.4) and (3.31,-0.3) .. (0,0) .. controls (3.31,0.3) and (6.95,1.4) .. (10.93,3.29)   ;
		\draw   (513.4,279.8) .. controls (513.4,271.07) and (520.47,264) .. (529.2,264) .. controls (537.93,264) and (545,271.07) .. (545,279.8) .. controls (545,288.53) and (537.93,295.6) .. (529.2,295.6) .. controls (520.47,295.6) and (513.4,288.53) .. (513.4,279.8) -- cycle ;
		\draw    (509,267.6) -- (454.75,237.57) ;
		\draw [shift={(453,236.6)}, rotate = 28.97] [color={rgb, 255:red, 0; green, 0; blue, 0 }  ][line width=0.75]    (10.93,-3.29) .. controls (6.95,-1.4) and (3.31,-0.3) .. (0,0) .. controls (3.31,0.3) and (6.95,1.4) .. (10.93,3.29)   ;
		\draw  [dash pattern={on 4.5pt off 4.5pt}]  (313,315.6) -- (209,317.56) ;
		\draw [shift={(207,317.6)}, rotate = 358.92] [color={rgb, 255:red, 0; green, 0; blue, 0 }  ][line width=0.75]    (10.93,-3.29) .. controls (6.95,-1.4) and (3.31,-0.3) .. (0,0) .. controls (3.31,0.3) and (6.95,1.4) .. (10.93,3.29)   ;
		\draw   (170.4,317.8) .. controls (170.4,310.18) and (176.58,304) .. (184.2,304) .. controls (191.82,304) and (198,310.18) .. (198,317.8) .. controls (198,325.42) and (191.82,331.6) .. (184.2,331.6) .. controls (176.58,331.6) and (170.4,325.42) .. (170.4,317.8) -- cycle ;
		\draw    (169,305.6) -- (145.04,266.31) ;
		\draw [shift={(144,264.6)}, rotate = 58.63] [color={rgb, 255:red, 0; green, 0; blue, 0 }  ][line width=0.75]    (10.93,-3.29) .. controls (6.95,-1.4) and (3.31,-0.3) .. (0,0) .. controls (3.31,0.3) and (6.95,1.4) .. (10.93,3.29)   ;
		\draw   (95,223) -- (155,223) -- (155,256.6) -- (95,256.6) -- cycle ;
		\draw    (130,218.6) -- (106.53,133.53) ;
		\draw [shift={(106,131.6)}, rotate = 74.58] [color={rgb, 255:red, 0; green, 0; blue, 0 }  ][line width=0.75]    (10.93,-3.29) .. controls (6.95,-1.4) and (3.31,-0.3) .. (0,0) .. controls (3.31,0.3) and (6.95,1.4) .. (10.93,3.29)   ;
		\draw   (217,249.8) .. controls (217,239.42) and (232.45,231) .. (251.5,231) .. controls (270.55,231) and (286,239.42) .. (286,249.8) .. controls (286,260.18) and (270.55,268.6) .. (251.5,268.6) .. controls (232.45,268.6) and (217,260.18) .. (217,249.8) -- cycle ;
		\draw  [dash pattern={on 0.84pt off 2.51pt}]  (298,257) -- (376.02,269.29) ;
		\draw [shift={(378,269.6)}, rotate = 188.95] [color={rgb, 255:red, 0; green, 0; blue, 0 }  ][line width=0.75]    (10.93,-3.29) .. controls (6.95,-1.4) and (3.31,-0.3) .. (0,0) .. controls (3.31,0.3) and (6.95,1.4) .. (10.93,3.29)   ;
		\draw  [dash pattern={on 0.84pt off 2.51pt}]  (220,237.6) -- (129.7,181.16) ;
		\draw [shift={(128,180.1)}, rotate = 32.01] [color={rgb, 255:red, 0; green, 0; blue, 0 }  ][line width=0.75]    (10.93,-3.29) .. controls (6.95,-1.4) and (3.31,-0.3) .. (0,0) .. controls (3.31,0.3) and (6.95,1.4) .. (10.93,3.29)   ;
		
		\draw (41,63.6) node [anchor=north west][inner sep=0.75pt]  [font=\footnotesize] [align=left] {helicopter with the \\shortest remaining\\rest time};
		\draw (73,22) node [anchor=north west][inner sep=0.75pt]  [font=\footnotesize] [align=left] {$\displaystyle t\ =\ 0$};
		\draw (236,76.6) node [anchor=north west][inner sep=0.75pt]   [align=left] {{\footnotesize closest rest}\\{\footnotesize base}};
		\draw (420,75) node [anchor=north west][inner sep=0.75pt]  [font=\footnotesize] [align=left] {does helicopter have \\time enough to perform \\a working circuit?};
		\draw (320,19) node [anchor=north west][inner sep=0.75pt]  [font=\footnotesize] [align=left] {no};
		\draw (380,153) node [anchor=north west][inner sep=0.75pt]  [font=\footnotesize] [align=left] {yes};
		\draw (340,201) node [anchor=north west][inner sep=0.75pt]  [font=\footnotesize] [align=left] {helicopter perfoms\\work circuit};
		\draw (336,291) node [anchor=north west][inner sep=0.75pt]  [font=\footnotesize] [align=left] {does helicopter have \\time enough to perform \\a working circuit?};
		\draw (520,272) node [anchor=north west][inner sep=0.75pt]  [font=\footnotesize] [align=left] {yes};
		\draw (177,310) node [anchor=north west][inner sep=0.75pt]  [font=\footnotesize] [align=left] {no};
		\draw (97,226) node [anchor=north west][inner sep=0.75pt]  [font=\footnotesize] [align=left] {go to rest\\base};
		\draw (225,240) node [anchor=north west][inner sep=0.75pt]  [font=\footnotesize] [align=left] {increase $\displaystyle t$};
			\end{tikzpicture}
	\label{fig:definicion_aeronaves}
\end{figure}
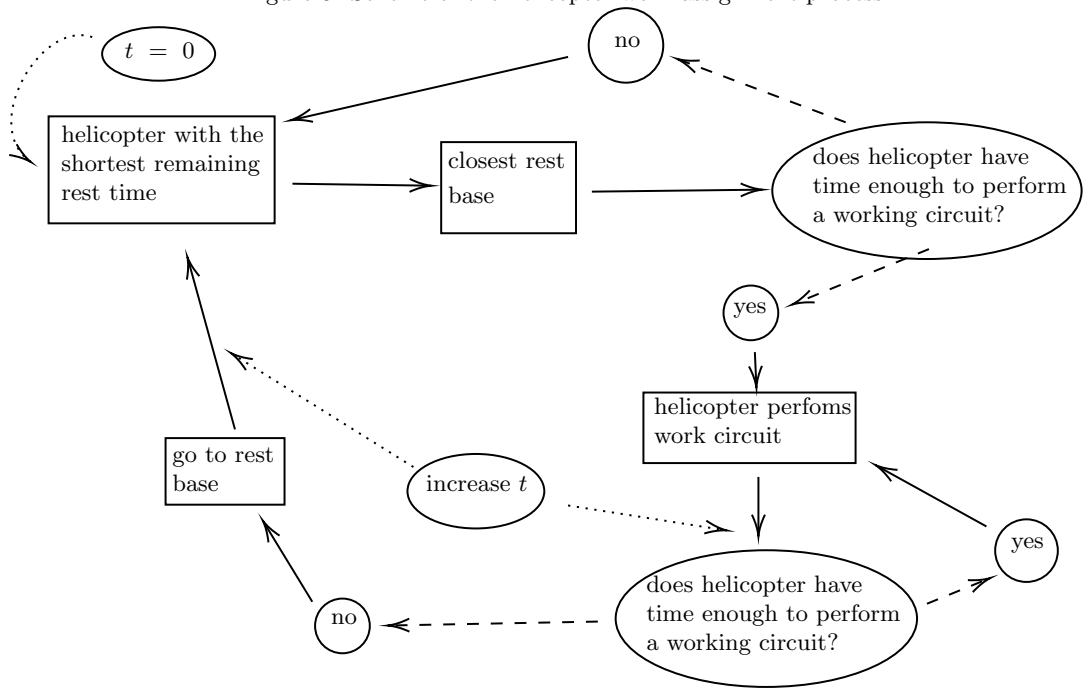

Figure \ref{fig:definicion_aeronaves} shows a scheme of the helicopter work assignment process.

 The helicopter that is ready to start working most quickly is chosen, taking into account that, in case of a tie, one of the better ones will be randomly selected. The rest base closest to the point of the wildfire corresponding to the time interval is selected, and it is evaluated if the helicopter can perform a work circuit, i.e., fly from the starting point to the water loading point, load water, fly to the wildfire, drop water and return to the base. The helicopter works at the water loading and dropping points assigned by the main trajectory to which it belongs as long as it has time available. When there is not time enough, it goes to the assigned base to rest and the next helicopter is selected.

When all the helicopters on the main trajectory have been selected for work, the first helicopter in the list is selected again for working once it has taken its mandatory rest.

Once the two stages have been applied, the feasibility of the solution is checked. In case of infeasibility, several movements, which depend on the violated constraints, can be applied. The movements are applied until a feasible solution is found. Table \ref{tab:feasibility_movements} describes all the movements implemented to ensure the feasibility of the solution.

\begin{table}[h]
	\centering
	\footnotesize
	\begin{tabularx}{\textwidth}{|c|X|X|}
		\hline
		\textbf{Id} & \textbf{Infeasibility}  & \textbf{Movement description} \\
		\hline
		MF1 & Capacity of rest bases exceeded.  & Modify the work schedule of a helicopter which rests in a base whose capacity is exceeded to assign it to other base. \\
		\hline
		MF2 & Water capacity of a water point is exceeded or helicopters capacity of a water point is exceeded. & Select a main trajectory which has been assigned the water load whose capacity is exceeded and replace the water load point. Update the work of all helicopters belonging to this main trajectory.\\
		\hline
		MF3 & Overlapping between helicopters on the same main trajectory. & Select one of the coinciding helicopters and anticipate or delay its working period.\\
		\hline
	\end{tabularx}
	\caption{Movements to obtain solutions' feasibility.}
	\label{tab:feasibility_movements}
\end{table}

\subsection{Solution improvement} \label{sec:improvement}
Once the initial feasible solution has been obtained, the local search phase begins. In this phase, solutions in the neighborhood of the current solution are explored. To this aim, the movements described in Table \ref{tab:optimal_movements} will be applied. These movements can be divided into two groups: the ones related to main trajectories and the movements associated with the helicopters. 

\begin{table}[h]
	\centering
	\footnotesize
	\begin{tabularx}{\textwidth}{|c|X|X|}
		\hline
		\textbf{Id} & \textbf{Movement}  & \textbf{Movement description} \\
		\hline
		MI1 & Change the water loading/dropping point.  & Assign a new water loading/dropping point to a main trajectory in a selected period. Update the work schedule of every helicopter belonging to the main trajectory. \\
		\hline
		MI2 & Change the water loading and dropping points. & Assign a new water dropping point to a main trajectory in a selected period. Select the water loading point closest to the newly selected water dropping point and assign it to the main trajectory in that period. Update the work schedule of each helicopter belonging to the main trajectory.\\
		\hline
		MI3 &Helicopter starts working. & Select a helicopter that is not in the current solution and force it to start working.\\
		\hline
		MI4 &Helicopter goes to rest. & Remove a helicopter from the current schedule. \\
		\hline
		MI5 &Helicopter backs to work. & Choose a helicopter which takes a rest longer than the minimum set and force it to go back to work when it finishes the minimum mandatory rest. \\
		\hline
		MI6 &Helicopter changes its working period. & Select a helicopter and change the time at which it starts to work. \\
		\hline
		MI7 &Helicopter changes the duration of a working circuit. & Select a helicopter and, each time it completes a water load-drop circuit, choose with some probability whether to return the helicopter to rest or continue working. \\   
		\hline
		MI8 &Extending the rest period of a helicopter. & Select a helicopter which is resting and extend its rest period to a random interval. \\
		\hline
	\end{tabularx}
	\caption{Movements to improve the solutions.}
	\label{tab:optimal_movements} 
\end{table}

Due to the complexity of the problem, it is necessary to check if the generated solutions are feasible. In case of infeasibility, movements of Table \ref{tab:feasibility_movements} will be considered until feasibility is reached. If it is not possible to restore feasibility, the movement is discarded.

\begin{figure}[h!]
	\caption{Example of the application of movements in Table \ref{tab:optimal_movements}.}                          
	\centering
	\includegraphics[width=0.8\textwidth]{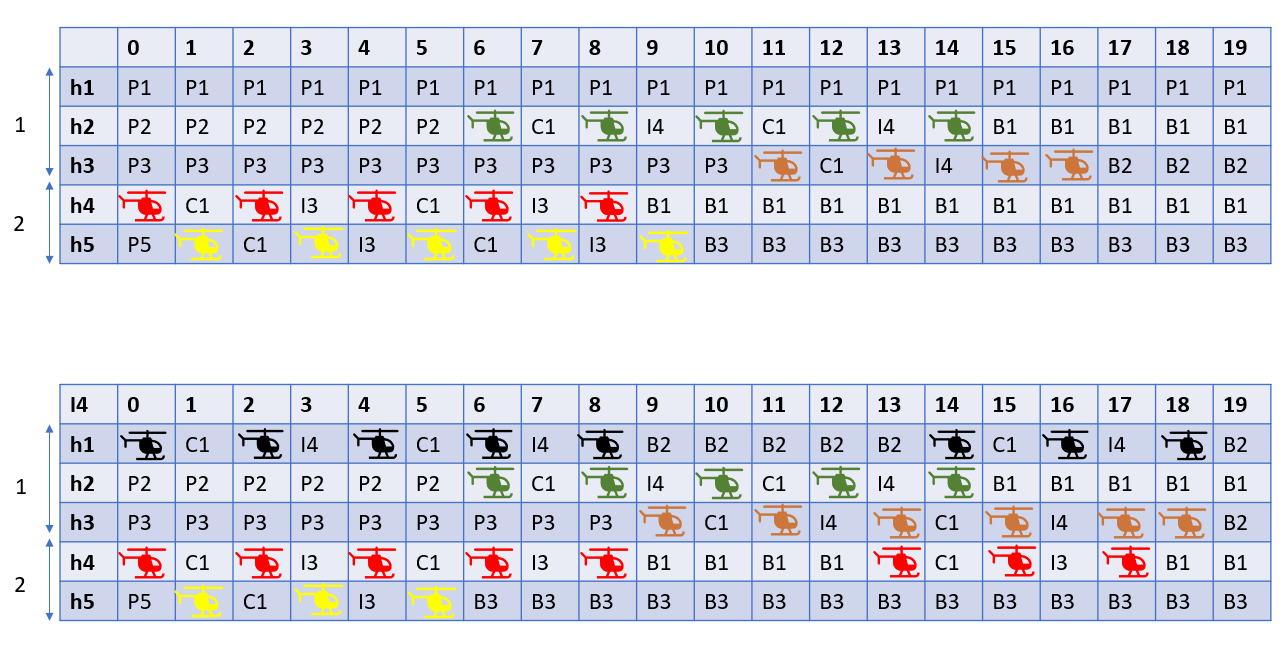}
	
	\label{fig:mov_mejora}
\end{figure}

An example of the application of the movements described in Table \ref{tab:optimal_movements} is shown in Figure \ref{fig:mov_mejora}. Each row of each table in Figure \ref{fig:mov_mejora} represents the work of a helicopter, and the columns represent time intervals. The letter P means that the helicopters are at their initial position, whereas the letters I, C and B represent the wildfire points, the water load points and the rest bases, respectively. The helicopter icon indicates that it is flying between positions. There are two main trajectories: one of them is composed of helicopters h1, h2 and h3 and the other takes into account h4 and h5.

The top of Figure \ref{fig:mov_mejora} represents the current solution, which is feasible. For simplicity, it is assumed that in this time period the wildfire nodes I1 and I2 have an associated efficiency of 5 and nodes I3 and I4 have an associated efficiency of 10. In the bottom of Figure \ref{fig:mov_mejora}, it can be seen the solution obtained after the implementation of several movements of those described in Table \ref{tab:optimal_movements}. Hence,  movement MI3 is applied to helicopter h1, which starts working at time interval t=0. The movement MI6 is applied to helicopter h3, causing it to start work at t=9. Movement MI5 is applied to helicopter h4 and MI7 movement is applied to helicopter h5, causing it to rest at time interval t=6 instead of exhausting its consecutive working time.

At a glance it can be seen that the main objective value of the problem, related to the efficiency of the water drops made, increases after the movements. In fact, the number of water drops increases from 7 to 11 drops. Thus, let us assume that all aircraft have the same capacity, e.g. 1000 liters per helicopter. Under this assumption, the efficiency of water drops would change from 70000 liters to 110000 liters.

\subsection{The G+SA and G+ILS algorithms}
Since it is not always possible to obtain good solutions for realistic instances in reasonable times with the commercial solver, two algorithms based on metaheuristic techniques have been developed.

Both algorithms start with an initial solution obtained with the procedure described in Subsection \ref{sec:initial_solution} and consider in the improvement phase the movements described in Subsection \ref{sec:improvement}. However, they differ in the metaheuristic techniques employed in the improvement phase. Whereas the G+SA algorithm combines the Greedy procedure to get the initial solution with an improvement phase based on the Simulated Annealing procedure, the G+ILS algorithm employs the Iterated Local Search metaheuristic to improve the initial solution.

As indicated in the pseudocode \ref{alg:algorithm_GSA}, the improvement phase of the G+SA algorithm is based on the cooling process defined in the classical Simulated Annealing technique (\cite{handbook_met}, \cite{kirkpatrick}, \cite{kirkpatrick2}), with a minor modification. To prevent from getting stuck at a local optimum, if the current solution is not improved after a number of iterations, a new initial solution is generated. Moreover, according to what was proposed in \cite{eglsese}, during the execution of the metaheuristic, the best result obtained is always stored.

In pseudocode \ref{alg:algorithm_GSA}, $x_0$ is current solution. In the first step, in order to define a new $x_0$, it is defined each main trajectory and then the helicopters work. Finally, it is checked that the solution is feasible, as explained in Subsection \ref{sec:initial_solution}. To obtain a new solution from $x_0$, it is applied a movement from Table \ref{tab:optimal_movements}, as explained in Subsection \ref{sec:improvement}. $x^*$ is the best solution, $\alpha$ is the cooling parameter, $T_k$ is the current temperature and $T_m$ is the final temperature. $L$ is
the limit of movements without improve the current solution, $LB$ is the limit of
movements without improve the best solution, $L_k$ is the counter of movements
without improve the current solution and $LB_k$ is the counter of movements
without improve the best solution.

\begin{minipage}{0.46\textwidth}
	\begin{algorithm}[H]		
		\centering
		\caption{G+SA}\label{alg:algorithm_GSA}
		\begin{algorithmic}[1]
			\State {Define $x_0$} 
			\State {$x^* := x_0$}
			\While{$T_k > T_m$ \textbf{and} $L_k < L$}
			\State {Generate $x_1$ from $x_0$}
			\If{$f(x_1) > f(x_0)$}
			\State{ $x_0 := x_1$}
			\Else \State { with probability $e^{\frac{f(x_1)-f(x_0)}{T_k}}$:}
			\State {$x_0 := x_1$}
			\EndIf
			\If{$x_0 = x_1$}
			\State {$L_k := 0$}
			\If{$f(x_1) > f(x^*)$}
			\State {$x^* := x_1$}
			\State {$LB_k := 0$}
			\Else
			\State {$LB_k = LB_k + 1$}
			\EndIf
			\Else
			\State {$L_k = L_k + 1$}
			\State {$LB_k = LB_k + 1$}
			\EndIf
			\State {$T_k :=\alpha T_k $}
			\If{$LB_k \geq LB$}
			\State {Define $x_0$} 
			\EndIf
			\EndWhile
		\end{algorithmic}
	\end{algorithm}
\end{minipage}
\hfill
\begin{minipage}{0.46\textwidth}
	\begin{algorithm}[H]
		\centering
		\caption{G+ILS}\label{alg:algorithm_GILS}
		\begin{algorithmic}[1]
			\State {Define $x_0$} 
			\While{$n < N_1$ \textbf{and} $L_k < L$}
			\State {Generate $x_1$ from $x_0$}
			\While{$n < N$}
			\State {Generate $x_2$ from $x_1$}
			\If{$f(x_2) > f(x_1)$}
			\State $x_1 := x_2$
			\State $f(x_1) := f(x_2)$
			\EndIf
			\State $n:= n+1$
			\EndWhile
			\If{ $f(x_2) > f(x_0)$}
			\State $x_0 := x_2$
			\State $f(x_0) := f(x_2)$
			\State $L_k := 0$
			\Else
			\State $L_k := L_k + 1$
			\EndIf
			\State $n := n+1 $
			\EndWhile
		\end{algorithmic}
	\end{algorithm}
\end{minipage}

\vspace*{0.3cm}
On the other hand, the improvement phase of G+ILS algorithm \ref{alg:algorithm_GILS} is related to the Iterated Local Search scheme (\cite{handbook_met2}) with Hill Climbing.

\section{Results}
This section deals with the analysis of the results obtained with both the commercial solver Gurobi and the algorithms based on the metaheuristic techniques (G+SA and G+ILS algorithms).

The good performance of these tools was tested with instances obtained from real data. Hence, real data of different models of helicopters have been used and their speeds were determined according to the average cruise speeds in the historical data of flights or from their technical data sheet. The water load points were provided by the Galician government \textit{teXunta de Galicia}, as well as the list of rest bases, which are available in PLADIGA (\cite{Pladiga}). Finally, wildfire data from two real wildfires, one small and one larger, were considered to get the coordinates on the wildfire evolution.

Thus, the combination of these inputs result in several data instances, which were divided into small, medium and big instances. Although the small instances are the less realistic, they are also useful since in these data Gurobi achieves the optimal solution. This is not the case of the medium and big instances (the most realistic ones), where Gurobi is not able to guarantee the optimality or even able to find a feasible solution. Since both the G+SA and G+ILS algorithms are capable of finding feasible solutions in reasonable times for all the instances, we will analyze if they can reach the optimal solution in the small instances and we will study their behavior in the larger ones.

The computational study with Gurobi has been made on a computer with a  processor Intel Xeon Gold 6126 12-Core, 384GB RAM. As far as the G+SA and G+ILS algorithms is concerned, all runs have been made on a computer with a processor 11th Gen Intel Core i7-1165G7 @2.80 GHz, 16GB RAM.

\subsection{Small data instances}
Small instances consist of 10 data files. As it can be seen in Table \ref{tab:computational_small}, the number of helicopters varies from 2 to 7, according to this, the number of main trajectories varies from 1 to 3. In all the files, there are 3 possible water dropping points, 5 water loading points and 5 rest bases. The time schedules varies from 2 to 4 hours, with time intervals of 2, 5 and 10 minutes. Table \ref{tab:computational_small} shows the description of these data files (number of helicopters, number of main trajectories, hours of scheduling, duration in minutes of the intervals into which the scheduling is divided) and the computational times (in minutes) obtained with the different methods.

\begin{table}[h]
	\footnotesize
	\centering	
	\begin{tabularx}{\textwidth}{|r|X|X|X|X|X|X|X|}
		\hline
		& & & & & \multicolumn{3}{c|}{Computational times} \\
		\hline
		 Set & Helicopters & Main trajectories & Hours & Time intervals &  Gurobi & G+SA & G+ILS \\
		\hline
		S1 & 2 & 2 & 2 & 5 & 0.030 & 0.005 & 0.003 \\

		S2 & 3 & 2 & 2 & 5 & 0.055 & 0.131 & 0.046 \\

		S3 & 5 & 3 & 2 & 5 & 0.095 & 0.722 & 0.471 \\

		S4 & 7 & 3 & 2 & 5 & 0.164 & 0.924 & 0.204 \\

		S5 & 3 & 2 & 4 & 10 & 8.660 & 0.110 & 0.098 \\

		S6 & 2 & 2 & 4 & 5 & 23.014 & 1.374 & 0.171 \\

		S7 & 3 & 2 & 2 & 2 & 52.107 & 1.200 & 1.179 \\

		S8 & 2 & 1 & 2 & 2 & 61.496 & 0.472 & 0.342\\

		S9  &5 & 3 & 2 & 2 & 99.634 & 3.566 & 0.927 \\
 
		S10 & 7 & 3 & 4 & 10 & 580.03 & 1.079 & 0.556 \\
		\hline		
	\end{tabularx}
	\caption{Computational times of small data instances.}
	\label{tab:computational_small}
\end{table}

\begin{table}[h]
	\footnotesize
	\centering	
	\begin{tabularx}{\textwidth}{|r|X|X|X|X|}
		\hline
		Set & Number of Drops & Number of Flights & Number of Changes & Blank times  \\
		\hline
		S1 &  8 & 27 & 0 & 0  \\
		S2 & 12 & 40 & 0 & 0  \\
		S3 & 19 & 68 & 0 & 0 \\
		S4 & 27 & 90 & 0 & 0 \\
		S5 & 9 & 27 & 0 & 0 \\
		S6 & 15 & 43 & 0 & 0 \\
		S7 & 25 & 95 & 0 & 0 \\
		S8 & 14 & 59 & 0 & 0 \\
		S9  & 45 & 168 & 0 & 0 \\
		S10 & 21 & 68 & 0 & 0 \\		
 		\hline
	\end{tabularx}
	\caption{Results of small data instances.}
	\label{tab:small_result}
\end{table}

From Table \ref{tab:computational_small} we can see that when the data instances are very small all methods reach the optimum in times of less than 1 minute. However, when they start to increase, as is the case of datasets S6, S7, S8, S9 and S10 Gurobi takes too long, considering the time it is reasonable to expect a result in this context. On the other hand, both algorithms (G+SA and G+ILS) take less than 5 minutes in all cases. It can be observed that the G+ILS, in particular, takes less than 2 minutes in all cases tested. Although both algorithms perform well on these datasets, G+ILS shows a better result, reaching the optimal value in less computational time.

Table \ref{tab:small_result} shows the objective function terms of the same datasets in the optimal result (number of total aircraft drops, total number of intervals they are flying, number of changes in main trajectories and number of time intervals where no aircraft are working on the wildfire). In none of the calculated schedules there are blank times, in which no aircraft is working on the wildfire, not even in those where only 2 helicopters work (``Blank times'' column).  As these are small wildfires, in none of these solutions there are changes in the main trajectories followed by the helicopters.

Table \ref{tab:small_result} also shows how the duration of the intervals affects the quality of the solution. For example, comparing files S1 and S8, in which there are 2 helicopters working for 2 hours, it can be seen that in S8 the total number of drops is 14, while in S1 it is 8. This is due to the length of the intervals. In S8 it is 2 minutes, which allows to better model the real behavior of the helicopters, while with 5 minutes intervals the accuracy is lost.

\subsection{Medium data instances}
Medium data instances are more realistic than small ones. In these datasets it has not been possible to guarantee optimality with Gurobi, although Gurobi has provided feasible solutions  to compare them with the solutions of the G+SA and G+ILS. To this end, the optimization model was run with Gurobi for 24 hours, and the best feasible solution obtained was selected.

These data instances are composed of 20 datasets. In these datasets, four types of wildfires were considered, depending on the number of nodes and the number of times that wildfire evolves during the scheduling: 5 nodes and only 1 evolution; 3 nodes and 5 evolutions; 5 nodes and 5 evolutions and 10 nodes and 3 evolutions.

Instances with 5, 10, 15, 20 and 25 helicopters have been created to use realistic scenarios, as it is considered that with less than 5 aircraft it does not make sense to use a helicopter planning tool.

Due to the high number of load points and rest bases existing in Galicia, and based on the historical data, the number of water load points and rest bases was fixed, assuming that 10 water load points and 5 bases are enough for all wildfires.

After the executions, it was found that as the number of time interval increases, the problem becomes too large and Gurobi has memory problems and is unable, within 24 hours at most, to find feasible solutions. For this reason, all medium data instances correspond to schedules of 8 hours and 10 minutes time intervals.

\begin{table}[h]
	\footnotesize
	\centering	
	\begin{tabularx}{\textwidth}{|r|X|X|X|X|X|}
		\hline
		Set & Helicopters & Main trajectories & Drop points  &  MIPGap & Best Value \\
		\hline
		M1 & 5 & 3 & 10  & 0.582 &  18.878 \\
		\hline
		M2 & 10 & 3  & 10 & 10.605 & 4.058 \\
		\hline
		M3 & 15 & 5 & 10 & 1.435 & 9.299 \\
		\hline
		M4 & 20 & 5 & 10 & 2.796 & 4.741 \\
		\hline
		M5 & 25 & 5 & 10 & 4.462 & 5.241 \\
		\hline
		M6 & 5 & 3 & 18 & 3.587 &   7.972\\
		\hline
		M7 & 10 & 3 & 18 & 155.457  &3.748  \\
		\hline
		M8 & 15 & 3 & 18 &5.754  & 3.180 \\
		\hline
		M9 & 20 & 3 & 18 & 19.448 & 0.776 \\
		\hline
		M10 & 25 & 3 & 18 & 144.225 & -1.957\\
		\hline
		M11 & 5 & 3 & 30  & 1.561 &  2.943 \\
		\hline
		M12 & 10 & 3  & 30 & 7.680 & 1.492 \\
		\hline
		M13 & 15 & 5 & 30 & 15.372 & 0.724 \\
		\hline
		M14 & 20 & 5 & 30 & 19.093 & 0.435 \\
		\hline
		M15 & 25 & 5 & 30 & 15.581 & 0.455 \\
		\hline
		M16 & 5 & 3 & 40 & 5.895 & 0.581  \\
		\hline
		M17 & 10 & 3 & 40 & 11.698 & -0.907 \\
		\hline
		M18 & 15 & 3 & 40 & 87.833 &  0.032\\
		\hline
		M19 & 20 & 3 & 40 & 53.326 & 0.114 \\
		\hline
		M20 & 25 & 3 & 40 & 32.600 & 0.161 \\
		\hline		
	\end{tabularx}
	\caption{Results of medium data instances.}
	\label{tab:computational_medium}
\end{table}

Table \ref{tab:computational_medium} shows a description of medium data files (number of aircraft, number of main trajectories and number of wildfire nodes). It also shows the Relative MIP optimality GAP (MIPGap) and the best value of the objective function (Best Value) at the time when the execution is stopped.

\begin{figure}[h!]
	\caption{Results of Gurobi, G+SA and G+ILS in medium instances.}
	\centering
	\includegraphics[width=1\textwidth]{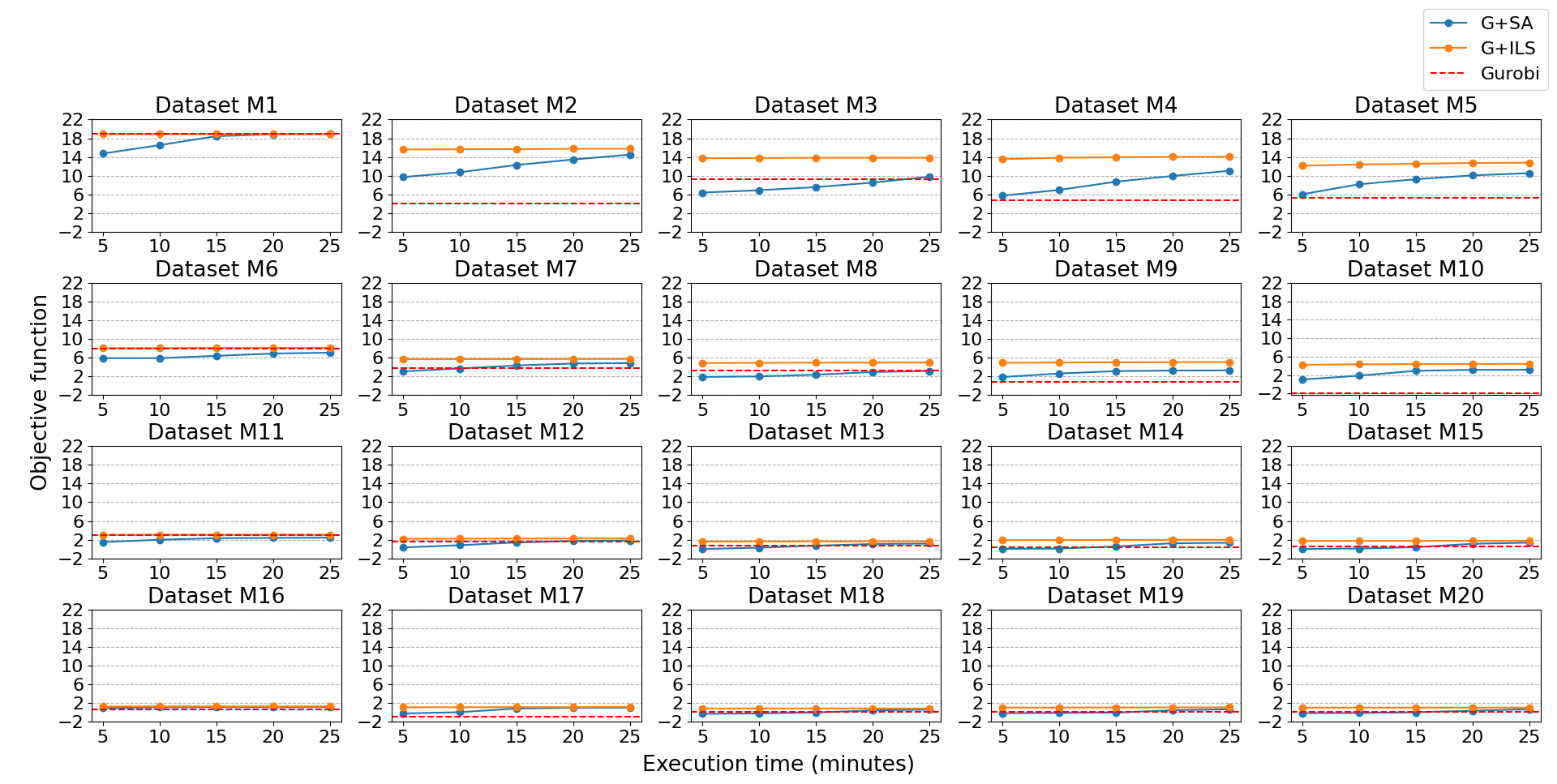}
	\label{fig:med_1}
\end{figure}

Figure \ref{fig:med_1} shows the results obtained with Gurobi (after a 24-hour run). The results obtained in Gurobi are compared with the results of G+SA and G+ILS after 5, 10, 15, 20 and 25 minutes, which are adequate times for a decision-support tool for the aerial coordinator. Both algorithms were run 10 times, obtaining the mean of the objective function.

Figure \ref{fig:med_1} exhibits the good performance of G+SA and G+ILS algorithms, which are 
far better than Gurobi. In fact, if we compare in terms of the objective function, the best solution obtained by Gurobi after 24 hours with the results of the algorithm after 25 minutes, G+ILS is clearly better than Gurobi for all files, whereas G+SA is superior in almost all files, excepting datasets M6 and M11.

On the other hand, G+ILS seems to be the best method by far, since it is able to find the better solution after 5 minutes of execution. Moreover, has a very stable behavior over the time.

To better illustrate the stability of the two algorithms, the Relative Percentage Difference (RDP) has been calculated. This measure indicates how close each method is to the best solution at each time interval tested. The RDP is defined as $ RDP = \frac{z' - z}{z' + \epsilon}$, where $z'$ is the best value obtained, $z$ is the value for the current solution, and $\epsilon$ is a negligible value greater enough to avoid problem computations when $z'$ is too close to zero.

\begin{figure}[h!]
	\caption{RDP values for both algorithms in medium instances.}
	\centering
	\includegraphics[width=1\textwidth]{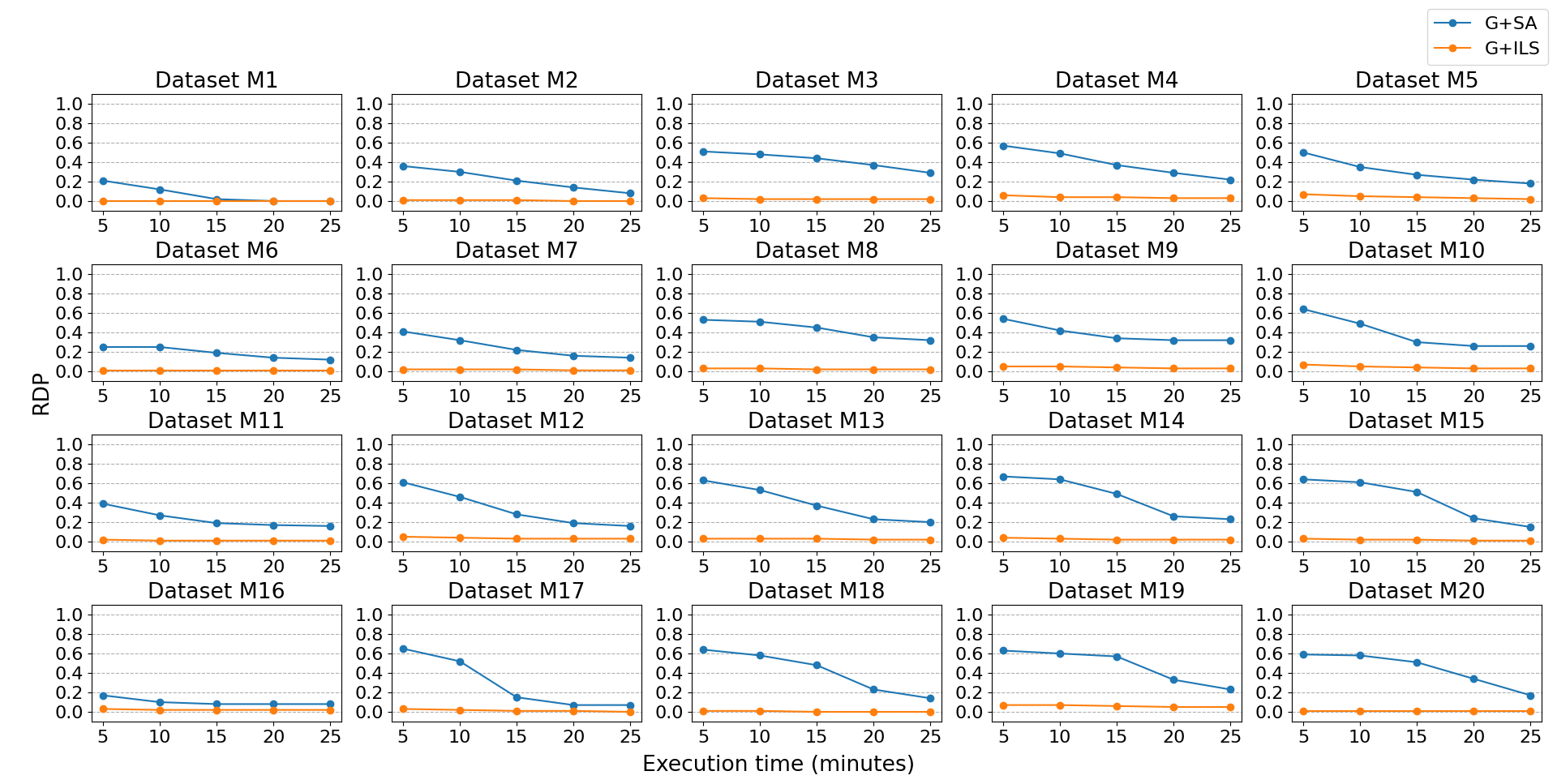}
	\label{fig:med_2}
\end{figure}

As in the previous case, $z$ represents the value of the average value of the objective function after 10 executions. A value of the RDP equal to zero indicates that the best solution has been found.
The results obtained for G+SA and G+ILS algorithms can be found in Figure \ref{fig:med_2}.

Figure \ref{fig:med_2} shows an unquestionable dominance of G+ILS over G+SA. Whereas G+ILS presents values close or equal to zero in all datasets at all execution times, G+SA only approaches zero in dataset M1 after 15 minutes. In addition, G+ILS is very stable over the time interval, while G+SA needs 20 to 25 minutes to achieve its best result.

\subsection{Big data instances}
Big data instances are the most realistic. It is necessary to use the developed tool  in daily schedules, up to 13 hours between sunrise and sunset. Thus, big data files are an extension of the input data considered in the medium instances, but the total time varies from 8 to 13 hours, with intervals of 1 and 5 minutes. Table \ref{tab:description_big} shows a description of all data files used. 

\begin{table}[h]
	\footnotesize
	\centering	
	\begin{tabularx}{\textwidth}{|r|X|X|X|X|X|}
		\hline
		Set & Helicopters & Main trajectories & Drop points  & Hours & Time intervals \\
		\hline
		B1 & 5 & 3 & 10  & 13 & 1 \\
		\hline
		B2 & 10 & 3  & 10 & 13 & 5 \\
		\hline
		B3 & 15 & 5 & 10 & 8 & 1 \\
		\hline
		B4 & 20 & 5 & 10 & 8 & 5 \\
		\hline
		B5 & 5 & 3 & 18 & 8 & 1 \\
		\hline
		B6 & 10 & 3 & 18 & 8 & 5 \\
		\hline
		B7 & 15 & 5 & 18 & 13  & 5 \\
		\hline
		B8  & 25 & 5 & 18 &  8 & 5 \\
		\hline
		B9 & 5 & 3 & 30 & 13 & 5 \\
		\hline
		B10 & 10 & 3 & 30 & 8 & 5 \\
		\hline
	\end{tabularx}
	\caption{Description of big data instances.}
	\label{tab:description_big}
\end{table}

\begin{figure}[h!]
	\caption{RDP values for both algorithms in big instances.}
	\centering
	\includegraphics[width=1\textwidth]{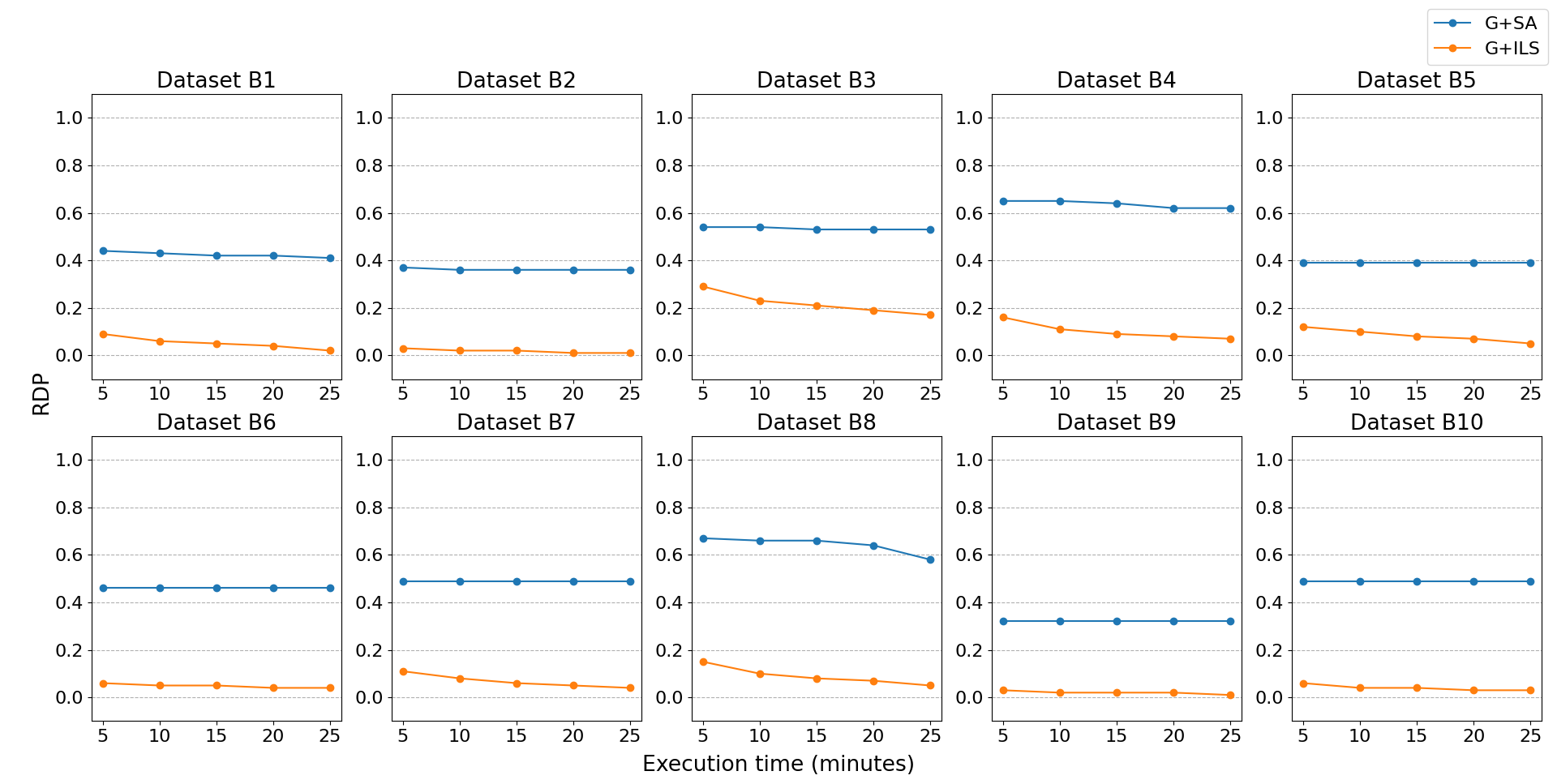}
	
	\label{fig:big_1}
\end{figure}

Since Gurobi has not been able to find feasible solutions for any of these instances after 24 hours, in this part only the results of G+SA and G+ILS will be analyzed.

Figure \ref{fig:big_1} shows the RDP values over the time intervals for G+SA and G+ILS algorithms. Again, greater stability is observed in G+ILS when compared to G+SA. In most of data files, the RDP of G+ILS is below of 0.2 after 5 minutes of execution. For instance B3 it takes slightly longer. Nevertheless, RDP value is below 0.4 in all datasets at 5 minutes of execution and, in the case of B3, is close to 0.2 after 10 minutes.

In \ref{appendix2} it is presented the result obtained with G+ILS for the data file B12.

\section{Conclusions}
This paper presents a MILP model to assist the wildfire coordinator in planning the work of firefighting helicopters in a large wildfire.

The developed model is close to the real performance of this type of aircraft in practice. It takes into account the fact that helicopters work according to main trajectories, which are defined by a water loading point, a water dropping point, as well as the number of helicopters that belong to the trajectory. The evolution of the wildfire over time is also taken into account, using the drop efficiency parameter, which helps the aerial coordinator to choose the areas where helicopters should drop water over the planning horizon. Helicopters are required to comply with the flight regulations in Spain, considering the capacity of the rest bases.

As real data are available, it has been possible to create realistic data instances, which allow the validation of the good performance of the methods.

First, the model was implemented in Gurobi solver, finding optimal solutions for small but not realistic data files. With the aim of providing alternative methods to find good solutions in reasonable times (less than 25 minutes) for realistic instances, they were developed two algorithms based on metaheuristic methods. Both algorithms have reached optimal solutions in small instances. But, although both methods have obtained feasible solutions in short periods of time for all data files and are competitive with Gurobi in small and medium instances, one of them is so far the better. The G+ILS algorithm, which combines a Greedy procedure to find the initial solution with the Iterated Local Search technique in the improvement phase, is able to find the better results after 5 minutes of execution in a wide variety of instances, with a very stable behavior.

It is worth mentioning that the methods provided take into account the characteristics and state of the helicopters at the time when they are executed. This enables them to run several times a day, which may be necessary if the evolution of the wildfire is not as expected, given its unpredictable nature.

A possible extension of the model would be to consider other types of extinguishing resources, such as aircraft or ground resources. Modeling the interaction of these resources would provide more information to the wildfire coordinator.






	\appendix 
	\section{}
	\label{appendix2}

	In this section it is presented a scheduling obtained with the ILS for the data instance B12. It is a 8-hours in 5-minutes schedule for 10 helicopters.

	Table \ref{tab:schedule_result} contains the time interval of each helicopter obtained in the computed scheduling.

	\begin{table}[h]
		\footnotesize
		\centering
		
		\begin{tabular}{|r|r|r|r|r|r|r|r|r|r|r|r|}
			\hline
			Helicopter    & \multicolumn{2}{c|}{W.I. 1}  & \multicolumn{2}{c|}{W.I. 2} & \multicolumn{2}{c|}{W.I. 3} & Number of drops \\
			\hline
			& Start & End & Start & End & Start & End &\\
			\hline
			h1 & 10:25 & 12:25 & 13:05 & 14:55 & 15:35 & 17:20 & 14\\
			\hline
			h2 & 10:55 & 12:55 & 13:35 & 15:20 & 16:00 & 17:45  & 15\\
			\hline
			h3 & 11:25 & 13:25 & 14:05 & 15:50 & 16:30 & 17:55 & 14 \\
			\hline
			h4 & 12:20 & 14:10 & 14:50 & 16:40 & 17:20 & 17:55 & 10 \\
			\hline
			h5 & 11:40 & 13:30 & 14:10 & 15:55 & 16:35 & 17:40 & 12\\
			\hline
			h6 & 11:30 & 13:25 & 14:05 & 15:50 & 16:30 & 17:55 & 13\\
			\hline
			h7 & 11:45 & 13:40 & 14:20 & 16:05 & 16:45 & 17:55 & 12\\
			\hline
			h8 & 11:20 & 13:15 & 13:55 & 15:40 & 16:20 & 17:45  &13\\
			\hline
			h9 & 10:20 & 12:15 & 12:55 & 14:40 & 15:20 & 17:05 & 14\\
			\hline
			h10 & 10:05 & 12:00 & 12:50 & 14:35 & 15:15 & 17:00 & 14  \\
			\hline       
		\end{tabular}
		\caption{Working times of each helicopter in the scheduled computed with the ILS.}
		\label{tab:schedule_result}
	\end{table}

	The value of the objective function of this result is equal to 9.7548, which is the result of:
	\begin{itemize}
		\itemsep0em 
		\item The sum of the efficiency of the water drops having into account the capacity of helicopters: 2162400.
		\item The number of time intervals in which helicopters are flying: 349.
		\item The number of time intervals in which helicopters are loading or dropping water: 262.
		\item The number of main trajectories changes: 11.
		\item The sum of all time intervals remaining for helicopters to work the maximum consecutive time: 38.
	\end{itemize}

\begin{figure}[h!]
		\caption{Helicopters' scheduling obtained with the Iterated Local Search for the Dataset B12.}
		\centering
		\includegraphics[width=1\textwidth]{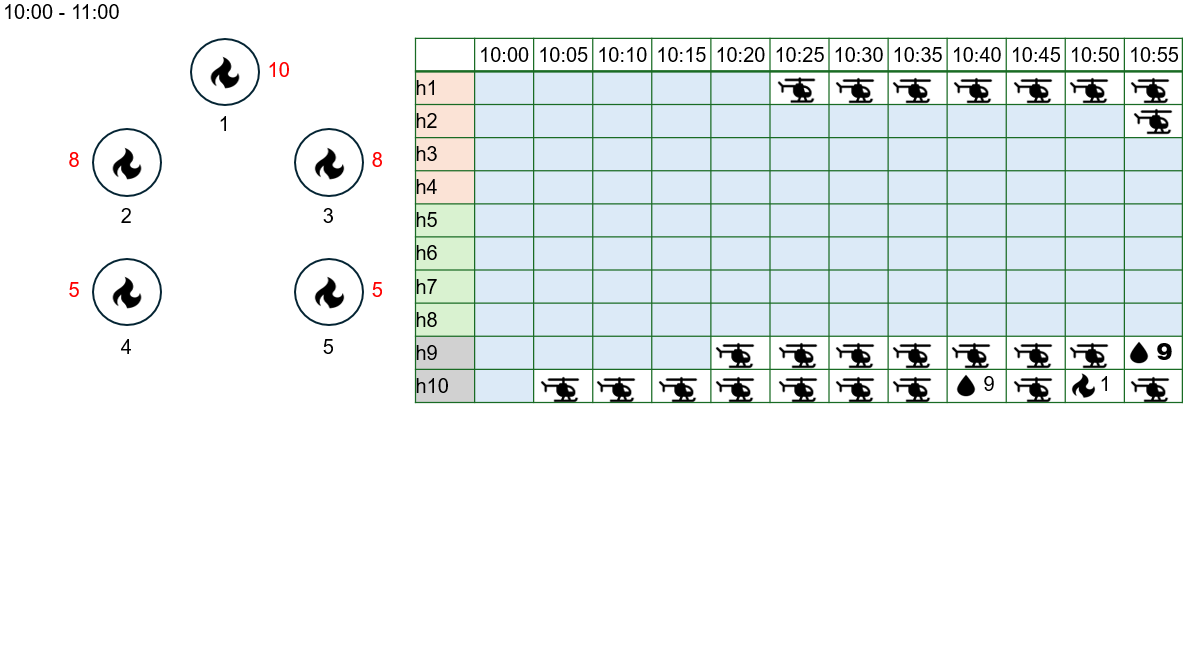}
		\includegraphics[width=1\textwidth]{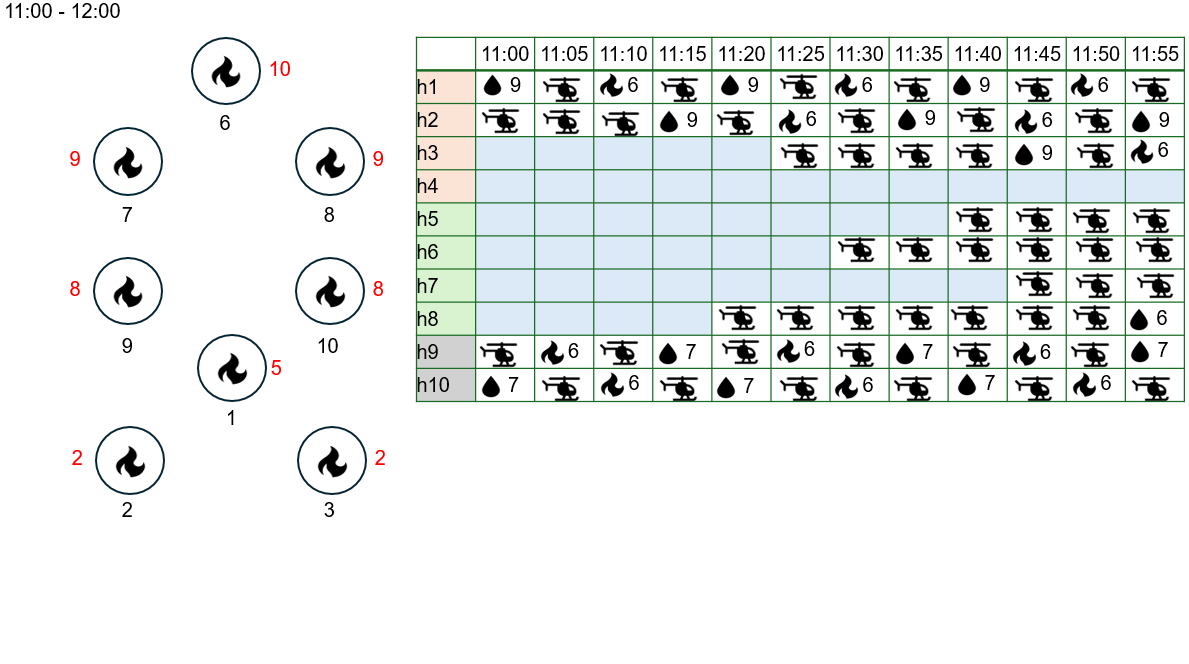}
		\label{fig:plan57_1}
	\end{figure}

	Figure \ref{fig:plan57_1} shows the complete helicopters' scheduling.
	
	The circles represent the nodes of the wildfire (the dropping point) with its efficiency associated (in red colour). The black number is an identifier of the node (there are 5 nodes with 5 evolution of the wildfire, so there are in total 30 nodes).
	
	The tables represent the scheduled obtained. Each row represents an helicopter, and the colour on the left represent the main trajectory at which their belong to (h1, h2, h3 and h4 belong to main trajectory 1; h5, h6, h7 and h8 belong to main trajectory 2; h9 and h10 belong to main trajectory 3).
	
	In the tables, the icon of an helicopter represents that the helicopter is flying from one position to another. The blue cells represent that the helicopter is in its initial position. The houses represent that the helicopter is in a rest base, the drops represent that it is loading water (the number is an indicator of the water loading point), the flame represent that the helicopter is dropping water in the wildfire (the number indicates in which point of the wildfire is dropping it).
	
	\begin{figure}[h!]
	\centering
	\includegraphics[width=1\textwidth]{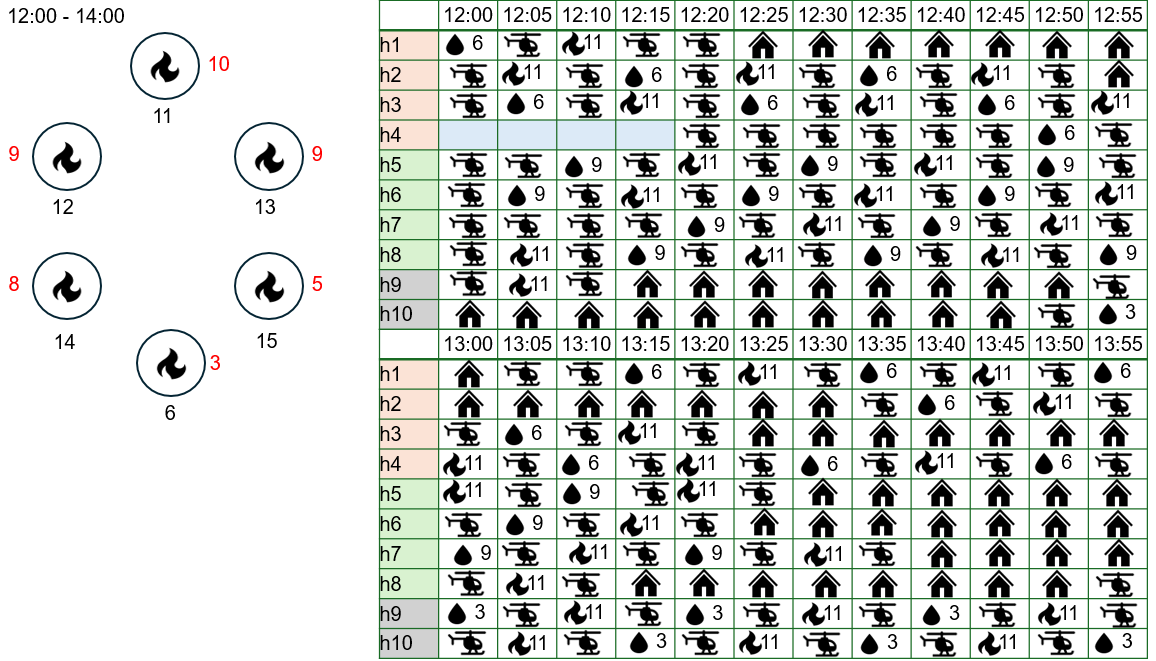}
	\includegraphics[width=1\textwidth]{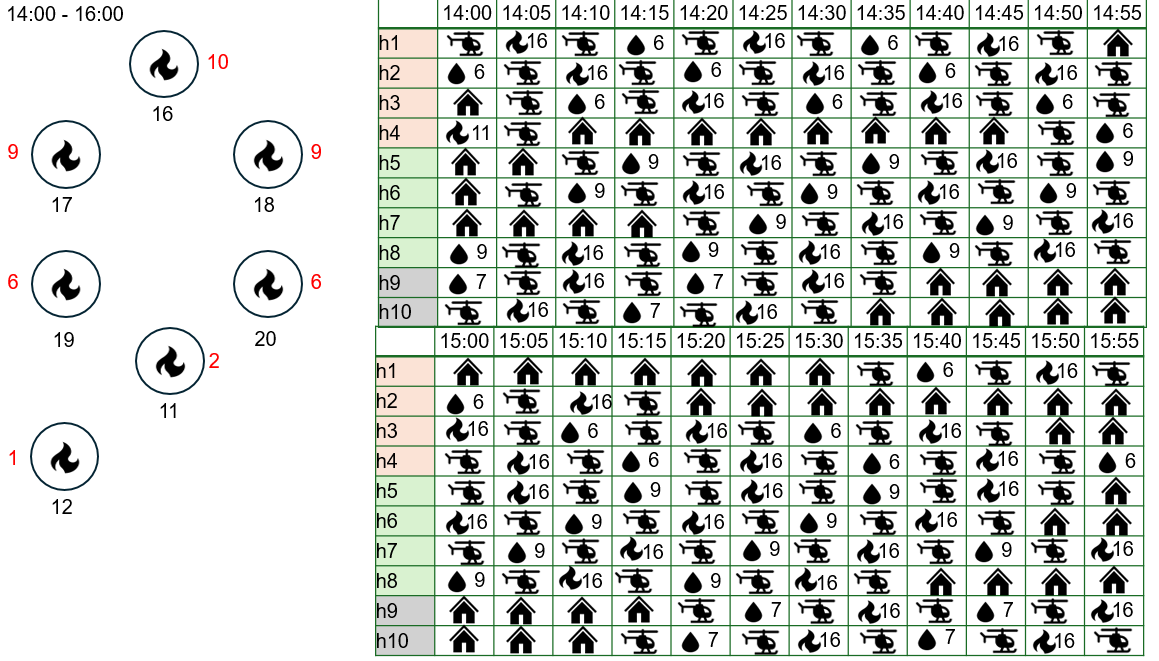}		
\end{figure}	
	\begin{figure}[h!]
		\centering
		\includegraphics[width=1\textwidth]{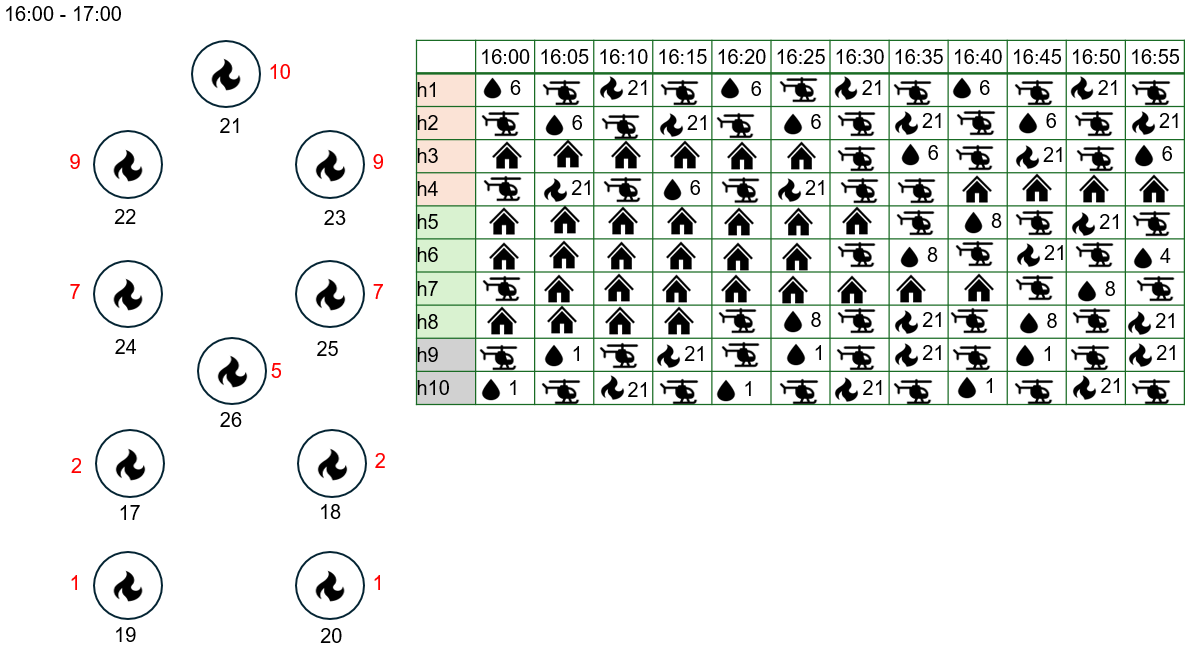}
		\\
		\includegraphics[width=1\textwidth]{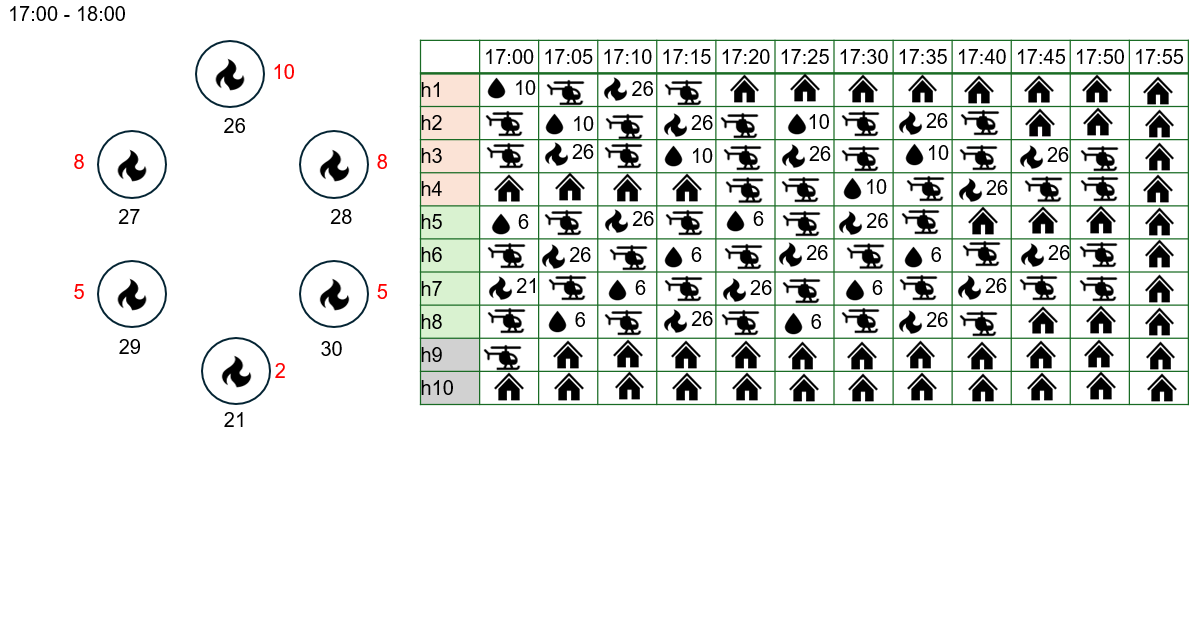}

\end{figure}

\end{document}